\documentclass[preprint,12pt]{elsarticle}

\usepackage[pdftex]{hyperref}
\biboptions{sort&compress,numbers,square}

\usepackage{amsfonts}
\usepackage{graphicx}
\usepackage{epsfig}
\usepackage{geometry}
\usepackage{tikz}
\usetikzlibrary{trees,positioning,shapes}
\usepackage{subcaption}
\usepackage{amsmath}
\usepackage{amsthm}

\usepackage{xcolor, colortbl}

\usepackage{algorithm}
\usepackage{algorithmicx, algpseudocode}

\newcommand{\rev}[1]{\textcolor{black}{#1}}
\newcommand{\inblue}[1]{\textcolor{blue}{#1}}
\newcommand{\inred}[1]{\textcolor{red}{#1}}

\usepackage{bm}
\usepackage{mathtools}

\usepackage{rotating}

\begin{document}

\begin{frontmatter}

\title{Isogeometric analysis with piece-wise constant test functions}

\author{Maciej Paszy\'{n}ski}

\address{Department of Computer Science, \\ AGH University of Science and Technology,
Krak\'{o}w, Poland \\
e-mail: maciej.paszynski@agh.edu.pl}

\begin{abstract}
We focus on the finite element method computations with higher-order $C^1$ continuity basis functions that preserve the partition of unity.
We show that the rows of the system of linear equations can be combined, and the test functions can be sum up to 1 using the partition of unity property at the quadrature points.
Thus, the test functions in higher continuity IGA can be set to piece-wise constants.
This formulation is equivalent to testing with piece-wise constant basis functions, with supports span over some parts of the domain. The resulting method is a Petrov-Galerkin formulation with piece-wise constant test functions.
This observation has the following consequences. The numerical integration cost can be reduced because we do not need to evaluate the test functions since they are equal to 1.
This observation is valid for any basis functions preserving the partition of unity property. It is independent of the problem dimension and geometry of the computational domain. It also can be used in time-dependent problems, e.g., in the explicit dynamics computations, where we can reduce the cost of generation of the right-hand side.
This summation of test functions can be performed for an arbitrary linear differential operator resulting from the Galerkin method applied to a PDE where we discretize with $C^1$ continuity basis functions.
The resulting method is equivalent to a linear combination of the collocations at points and with weights resulting from applied quadrature over the spans defined by supports of the piece-wise constant test functions.
\end{abstract}
	

\end{frontmatter}

\section{Introduction}

The main result of this paper can be summarized as follows.
We focus on finite element method discretization, with $C^1$ continuity basis functions, e.g. quadratic $C^1$ B-splines utilized in isogeometric analysis (IGA) \cite{iga1,iga2,iga3}.
Let us focus our attention on the one-dimensional Laplace equation for the simplicity of the presentation.
In this case, the Galerkin method involves the integrals $\int{B^x_{i,p}\Delta B^x_{j,p}}dx=-\int_{\Omega}{\nabla B^x_{i,p}\nabla B^x_{j,p}}dx$ (assuming zero boundary condition also for simplicity). If we use higher continuity basis functions, e.g., $C^1$ continuity B-splines, the approximation lives in a subspace of $H^2$, and if we integrate exactly, with proper numerical quadrature, these integrals are equal. In other words, the matrix of the system of linear equations resulting from the Galerkin method

\begin{eqnarray}
	\begin{bmatrix}
    -\int{{\nabla} B^x_{1,p}{\nabla}B^x_{1,p}}dx & \cdots &  -\int{{\nabla}B^x_{1,p}{\nabla}B^x_{N_x,p}}dx \\
    -\int{{\nabla}B^x_{2,p}{\nabla}B^x_{1,p}}dx & \cdots &  -\int{{\nabla}B^x_{2,p}{\nabla}B^x_{N_x,p}}dx\\   
    \vdots & \vdots & \vdots \\
-\int{{\nabla}B^x_{N_x,p}{\nabla}B^x_{1,p}}dx  
 & \cdots &  -\int{{\nabla}B^x_{N_x,p}{\nabla}B^x_{N_x,p}}dx \\
  \end{bmatrix} 
\begin{bmatrix}
  u_1 \\ u_2 \\ \vdots \\ u_{N_{x}}  \\
  \end{bmatrix} 
=    \begin{bmatrix}
\int {\cal RHS}(x) {B^x_1(x)} dx \\
\int {\cal RHS}(x) {B^x_2(x)} dx \\
\vdots \\
\int {\cal RHS}(x) {B^x_{N_x}(x)} dx \\
  \end{bmatrix}  \nonumber
\end{eqnarray}
with $C^1$ continuity basis functions have identical double precision values as the system not integrated by parts 
\begin{eqnarray}
	\begin{bmatrix}
    \int{B^x_{1,p}{\Delta}B^x_{1,p}}dx & \cdots &  \int{B^x_{1,p}{\Delta}B^x_{N_x,p}}dx \\
    \int{B^x_{2,p}{\Delta}B^x_{1,p}}dx & \cdots &  \int{B^x_{2,p}{\Delta}B^x_{N_x,p}}dx\\   
    \vdots & \vdots & \vdots \\
\int{B^x_{N_x,p}{\Delta}B^x_{1,p}}dx  
 & \cdots &  \int{B^x_{N_x,p}{\Delta}B^x_{N_x,p}}dx \\
  \end{bmatrix} 
\begin{bmatrix}
  u_1 \\ u_2 \\ \vdots \\ u_{N_{x}}  \\
  \end{bmatrix} 
=    \begin{bmatrix}
\int {\cal RHS}(x) {B^x_1(x)} dx \\
\int {\cal RHS}(x) {B^x_2(x)} dx \\
\vdots \\
\int {\cal RHS}(x) {B^x_{N_x}(x)} dx \\
  \end{bmatrix}  \nonumber
\end{eqnarray}
The matrices as well as the right-hand-sides of both systems are equal. The fluxes between elements are zero when we employ $C^1$ discretization. It does not matter which method we use for the generation of the system on the computer, and the resulting floating-point values will be the same (up to double precision round-off errors).

The second observation is that the system where we test the Laplace equation with $C^1$ B-splines can be transformed to the one where we have some piece-wise constant test functions ${\cal I_i}$, namely 
\begin{eqnarray}
	\begin{bmatrix}
    \int {\cal I}_{1}{\Delta}B^x_{1,p}dx & \cdots &  \int {\cal I}_{1}{{\Delta}B^x_{N_x,p}}dx \\
    \int {\cal I}_{2}{{\Delta}B^x_{1,p}}dx & \cdots &  \int {\cal I}_{2}{{\Delta}B^x_{N_x,p}}dx\\   
    \vdots & \vdots & \vdots \\
\int {\cal I}_{N_x}{{\Delta}B^x_{1,p}}dx  
 & \cdots &  \int {\cal I}_{N_x}{{\Delta}B^x_{N_x,p}}dx \\
  \end{bmatrix} 
\begin{bmatrix}
  u_1 \\ u_2 \\ \vdots \\ u_{N_{x}}  \\
  \end{bmatrix} 
=    \begin{bmatrix}
\int {\cal I}_{1}{{\cal RHS}(x)} dx \\
\int {\cal I}_{2}{{\cal RHS}(x)}dx \\
\vdots \\
\int {\cal I}_{N_x}{\cal RHS}(x)dx \\
  \end{bmatrix}. \nonumber
\end{eqnarray}
The details of the  derivation is described later in the paper. It is based on the idea of combining the rows of the matrix. The rows are combined in such a way that test functions sum up to 1, using the partition of unity property.

This observation has the following important consequences. First, the numerical integration cost will be reduced, since we do not need to integrate the test functions (we do not need to evaluate the test B-splines at quadrature points). Second, this observation does not depend on the selected quadrature points. Third, this transformation can be performed if we replace the Laplacian by any partial differential operator resulting from a PDE that can be solved with $H^2$ approximations with $C^1$ basis functions, preserving the partition of unity property. We selected B-splines for the simplicity of the presentation but there are several other options for discretization available \cite{TA,TB,TC,TD,HA,HB,HC,LA,LB,SA,SB,SC,SD}. The critical here is the partition of unity property. Fourth, this equivalence is also independent of the dimension of the problem, and it works in two or three-dimensions, or in the space-time formulations. Fifth, this equivalence is independent of the geometry of the computational domain and of the Jacobian of the transformation of the patch of elements into the master patch.
Recently, it is also possible to extend the $C^1$ continuity between patches of elements \cite{Patches}, so the equivalence with piece-wise constants also co be extended there.
Sixth, we end up with the integrals using the values of the trial functions at the quadrature points. 

It is like combining the collocation points \cite{col1,col2} at quadrature points with quadrature weights, over the spans of piece-wise constant test functions. 
The quadrature and the spans of the piece-wise constant test functions define the locations of the collocation points. Several collocation points are combined into one equation by the integration operator.

This observation speeds up also the explicit simulations with IGA since the integration of the right-hand side is cheaper.
The same logic applies to any basis functions that are globally $C^1$ and preserves the partition of unity property.

The structure of the paper is the following. We start in Section 2 from the one-dimensional derivation of the method. Next, we focus on the two-dimensional extension in Section 3. Finally, in Section 3, we illustrate the method with four numerical examples, the three-dimensional projection problem, the explicit dynamics simulation, the two-dimensional Laplace problem, and the isogeometric L2 projection of a bitmap. We summarize the paper in Section 4.

\section{One dimensional case}

Let us focus on the general PDE in the following form 
\begin{equation}
{\cal F}u={\cal RHS}
\end{equation}
defined over $\Omega=[a,b]$ interval. We partition the interval into $N_e$ finite elements.
Let us use the Galerkin method with $C^1$ continuity of the discretization.
We have the one dimensional B-spline basis functions
\begin{equation}
\{B^x_{i,p}(x)\}_{i=1,...,N_x}
\end{equation}
where $N_x=N_e+p$. We approximate the solution $u(x)\approx \sum_{i=1,...,N_x}u_i B^x_{i,p}(x)$.
We also test with B-splines.

If we have $C^1$ continuity of the trial basis functions and we use the exact quadrature during the integration, then the fact, if we integrate by parts or not, does not matter, the values in the matrix are the same, before or after the integration.
So let us focus on $C^1$ B-splines and test our PDE with B-splines, and we do not integrate by parts.

\begin{eqnarray}
	\begin{bmatrix}
    \int{B^x_{1,p}{\cal F}\left(B^x_{1,p}\right) }dx & \cdots &  \int{B^x_{1,p}{\cal F}\left( B^x_{N_x,p}\right)} dx \\
    \int{B^x_{2,p}{\cal F}\left(B^x_{1,p}\right)}dx & \cdots &  \int{B^x_{2,p}{\cal F}\left(B^x_{N_x,p}\right) }dx\\   
    \vdots & \vdots & \vdots \\
\int{B^x_{N_x,p}{\cal F}\left( B^x_{1,p}\right)} dx  
 & \cdots &  \int{B^x_{N_x,p}{\cal F}\left( B^x_{N_x,p}\right) }dx \\
  \end{bmatrix} 
\begin{bmatrix}
  u_1 \\ u_2 \\ \vdots \\ u_{N_{x}}  \\
  \end{bmatrix} 
=    \begin{bmatrix}
\int {\cal RHS}(x) {B^x_1} dx \\
\int {\cal RHS}(x) {B^x_2} dx \\
\vdots \\
\int {\cal RHS}(x) {B^x_{N_x}} dx \\
  \end{bmatrix}  \nonumber
\end{eqnarray}

Let us select any quadrature with points and weights $\{x_o,w_o\}_o$, resulting in the exact numerical integration.
At a given quadrature point $x_o$ we have $p+1$ non-zero B-spline functions.

We take our system of linear equations, and we replace the first row by the sum of rows $1,2,...,1+k$.
We also replace the second row by the sum of rows $1,2,...,2+k$.
Similarly, we replace row $r$ by the sum of rows $max(1,r-k),...,min(N_x,r+k)$ to the row $r=2,...,N_x-1$. 
Finally, we replace the last row by the sum of rows $N_x-k,...,N_x$. We get the system

\begin{eqnarray}
	\begin{bmatrix}
    \sum_{m=1,...,k+1}\int{B^x_{m,p}{\cal F}\left(B^x_{1,p}\right)}dx  & \cdots &  \sum_{m=1,...,k+1}\int{B^x_{m,p}{\cal F}\left(B^x_{N_x,p}\right)}dx \\
    \sum_{m=1,...,k+2}\int{B^x_{m,p}{\cal F}\left(B^x_{1,p}\right)}dx & \cdots &  \sum_{m=1,...,k+2}\int{B^x_{m,p}{\cal F}\left(B^x_{N_x,p}\right)}dx\\   
    \vdots & \vdots &  \vdots \\
\sum_{m=N_x-k,...,N_x}\int{B^x_{m,p}{\cal F}\left(B^x_{1,p}\right)}dx  
 & \cdots &  \sum_{m=N_x-k,...,N_x}\int{B^x_{m,p}{\cal F}\left(B^x_{N_x,p}\right)}dx \\
  \end{bmatrix} \begin{bmatrix}
  u_1 \\ u_2\\ \vdots \\ u_{N_{x}}  \\
  \end{bmatrix} \nonumber \\
=    \begin{bmatrix}
\sum_{m=1,...,k+1}\int {\cal RHS}(x) {B^x_{m,p}} dx \\
\sum_{m=1,...,k+2}\int {\cal RHS}(x) {B^x_{m,p}} dx \\
\vdots \\
\sum_{m=N_x-k,N_x}\int {\cal RHS}(x) {B^x_{m,p}} dx \\
  \end{bmatrix}  \nonumber
\end{eqnarray}

\begin{figure}
\includegraphics[scale=0.5]{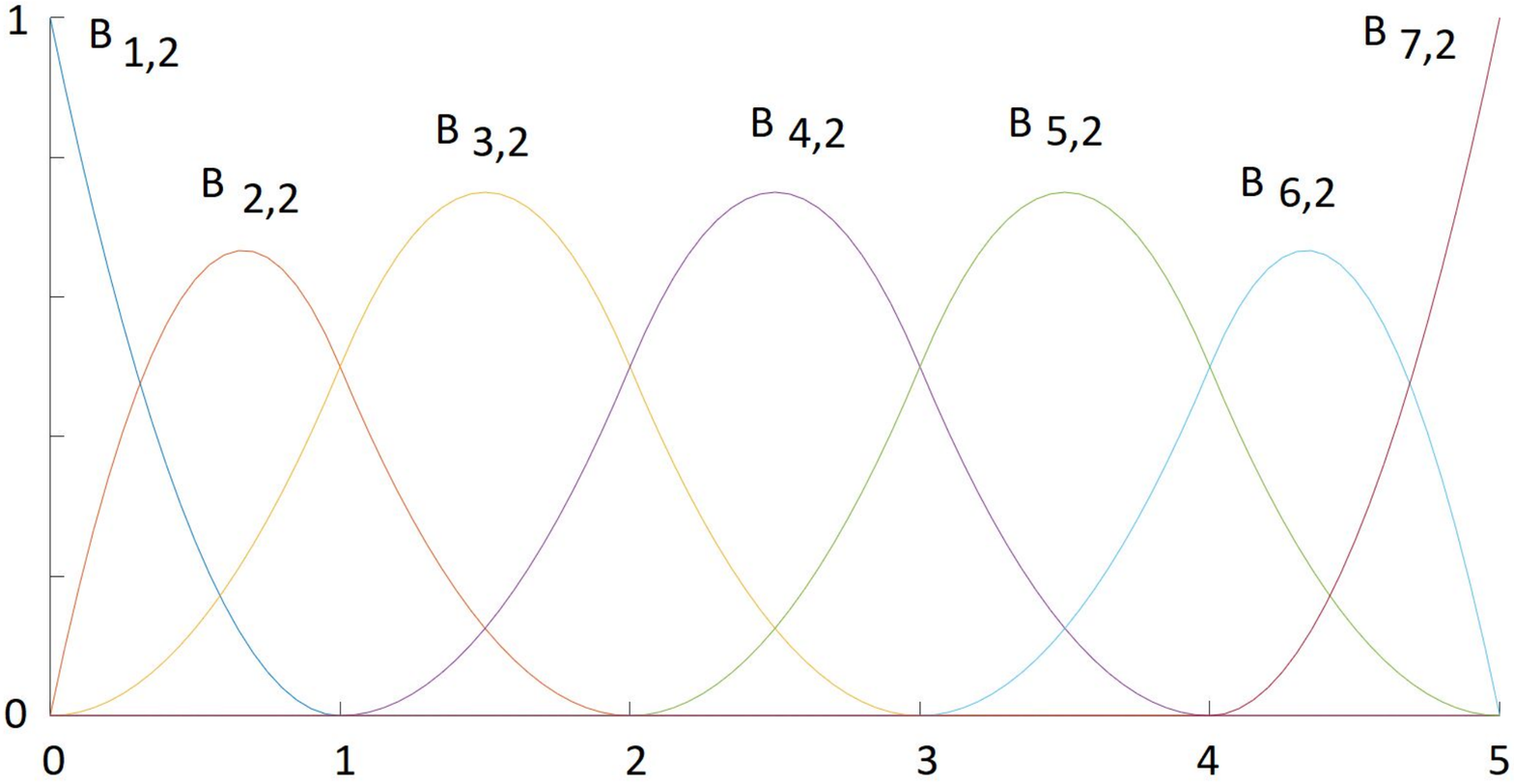}
\caption{B-splines span over [0 0 0 1 2 3 4 5 5 5] knot vector.}
\label{figure0}
\end{figure}

\begin{figure}
\includegraphics[scale=0.24]{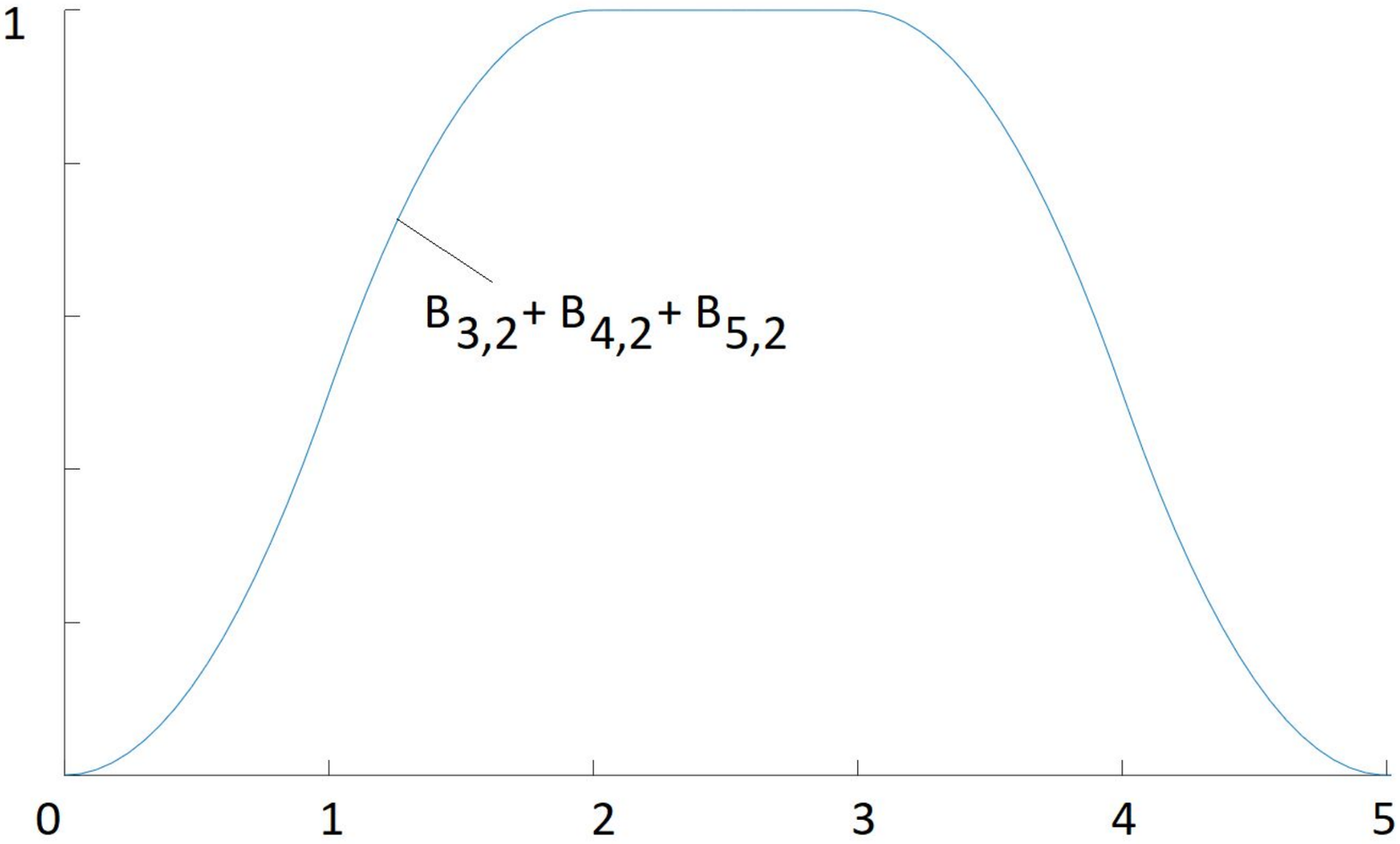}
\includegraphics[scale=0.24]{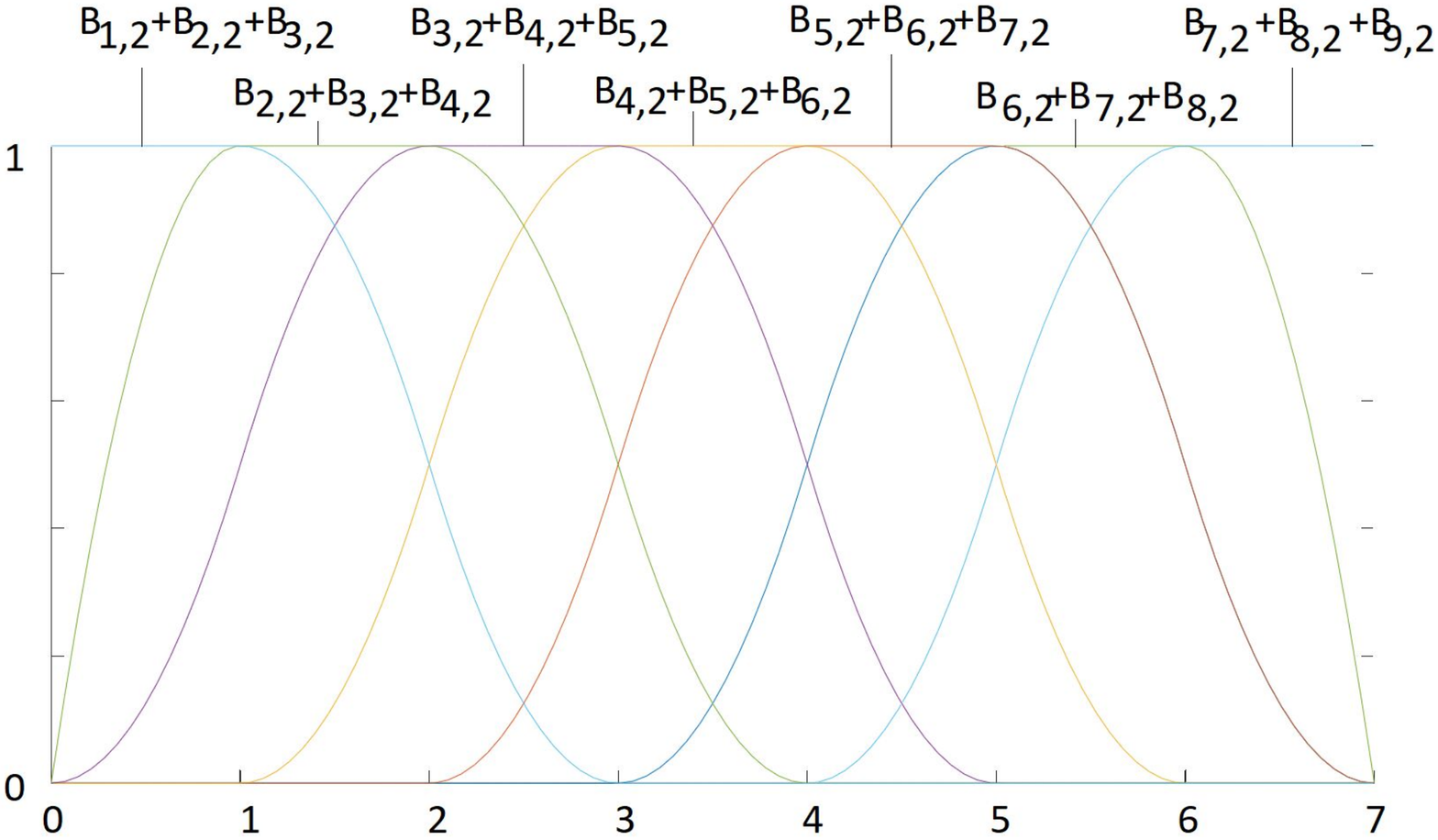}
\caption{One and all test functions obtained by summing up three consecutive B-splines.}
\label{fig:sum}
\end{figure}

Let us illustrate the matrix of the system by focusing on the following example.
Let us consider quadratic B-splines over 5 elements, defined by knot vector [0 0 0 1 2 3 4 5 5 5], which results in trial basis functions $B^x_{1,2},...,B^x_{7,2}$. Here, $B_1$ has support over [0,1], $B_2$ over [0,2], $B_3$ over [0,3], $B_4$ over [1,4], $B_5$ over [2,4], 
$B_6$ over [3,4], and $B_7$ over [4,5].  
We define now new test functions, by summing up three consecutive B-splines, 
$B^x_{i,2}+B^x_{i+1,2}+B^x_{i+2,2}$. The resulting new test functions are presented in Figure \ref{fig:sum}.

We can partition the integrals according to the supports of the basis functions.
We plot the entire matrix (in two blocks). For simplicity, we skip the superscript $^x$ in the notation.

{\footnotesize \begin{eqnarray}
\left[ \begin{array}{lll} 
\int_0^1(B_1+B_2){\cal F}B_1  &   \int_0^2 (B_1+B_2){\cal F}B_2 &   \int_0^3(B_1+B_2){\cal F}B_3 \\
\int_0^1(B_1+B_2+B_3){\cal F}B_1  &   \int_0^2 (B_1+B_2+B_3){\cal F}B_2 &   \int_0^3(B_1+B_2+B_3){\cal F}B_3 \\
\int_0^1(B_2+B_3+B_4){\cal F}B_1  &   \int_0^2 (B_2+B_3+B_4){\cal F}B_2 &   \int_0^3(B_2+B_3+B_4){\cal F}B_3 \\
\int_0^1(B_3+B_4+B_5){\cal F}B_1  &   \int_0^2 (B_3+B_4+B_5){\cal F}B_2 &   \int_0^3(B_3+B_4+B_5){\cal F}B_3 \\
\int_0^1(B_4+B_5+B_6){\cal F}B_1  &   \int_0^2 (B_4+B_5+B_6){\cal F}B_2 &   \int_0^3(B_4+B_5+B_6){\cal F}B_3 \\
\int_0^1(B_5+B_6+B_7){\cal F}B_1  &   \int_0^2 (B_5+B_6+B_7){\cal F}B_2 &   \int_0^3(B_5+B_6+B_7){\cal F}B_3 \\
\int_0^1(B_6+B_7){\cal F}B_1  &   \int_0^2 (B_6+B_7){\cal F}B_2 &   \int_0^3(B_6+B_7){\cal F}B_3 \\
	\end{array} \right. ... \nonumber
\end{eqnarray} 
\begin{eqnarray}
... \left. \begin{array}{llll} 
    \int_1^4(B_1+B_2){\cal F}B_4 &  \int_2^5(B_1+B_2){\cal F}B_5 &   \int_3^5(B_1+B_2){\cal F}B_6  & \int_4^5(B_1+B_2){\cal F}B_7 \\
    \int_1^4(B_1+B_2+B_3){\cal F}B_4 &  \int_2^5(B_1+B_2+B_3){\cal F}B_5 &   \int_3^5(B_1+B_2+B_3){\cal F}B_6 &   \int_4^5(B_1+B_2+B_3){\cal F}B_7 \\
    \int_1^4(B_2+B_3+B_4){\cal F}B_4 &  \int_2^5(B_2+B_3+B_4){\cal F}B_5 &   \int_2^5(B_2+B_3+B_4){\cal F}B_6  &   \int_4^5(B_2+B_3+B_4){\cal F}B_7 \\
    \int_1^4(B_3+B_4+B_5){\cal F}B_4 &  \int_2^5(B_3+B_4+B_5){\cal F}B_5 &   \int_2^5(B_3+B_4+B_5){\cal F}B_6 &  \int_2^5(B_3+B_4+B_5){\cal F}B_7 \\
    \int_1^4(B_4+B_5+B_6){\cal F}B_4 &  \int_2^5(B_4+B_5+B_6){\cal F}B_5 &   \int_2^5(B_4+B_5+B_6){\cal F}B_6  &   \int_2^5(B_4+B_5+B_6){\cal F}B_7 \\
    \int_1^4(B_5+B_6+B_7){\cal F}B_4 &  \int_2^5(B_5+B_6+B_7){\cal F}B_5 &   \int_2^5(B_5+B_6+B_7){\cal F}B_6  &   \int_2^5(B_5+B_6+B_7){\cal F}B_7 \\
    \int_1^4(B_6+B_7){\cal F}B_4 &  \int_2^5(B_6+B_7){\cal F}B_5 &   \int_2^5(B_6+B_7){\cal F}B_6  &   \int_2^5(B_6+B_7){\cal F}B_7 \\
	\end{array} \right] \nonumber
\end{eqnarray}}

We can organize these terms as follows

{\footnotesize{\begin{eqnarray}
\left[ \begin{array}{lll} 
\inred{\int_0^1(B_1+B_2){\cal F}B_1}  &   \inred{\int_0^2 (B_1+B_2){\cal F}B_2} &   \inred{\int_0^2(B_1+B_2){\cal F}B_3} \\
\int_0^1{\cal F}B_1  &  \int_0^1{\cal F}B_2+ \inred{\int_1^2 (B_2+B_3){\cal F}B_2} &  
\int_0^1{\cal F}B_3+ \inred{\int_1^2 (B_2+B_3){\cal F}B_3} + \inblue{\int_2^3B_3{\cal F}B_3} \\
\inred{\int_0^1(B_2+B_3){\cal F}B_1}  & \inred{\int_0^1(B_2+B_3){\cal F}B_2} +  \int_1^2 {\cal F}B_2 &   \inred{\int_0^1(B_2+B_3){\cal F}B_3} + \int_1^2 {\cal F}B_3 + \inred{ \int_2^3(B_3+B_4){\cal F}B_3} \\
\inblue{\int_0^1(B_3){\cal F}B_1}  &    \inblue{\int_0^1 (B_3){\cal F}B_2}+ \inred{\int_1^2 (B_3+B_4){\cal F}B_2} &   \inblue{\int_0^1(B_3){\cal F}B_3}+ \inred{\int_1^2(B_3+B_4){\cal F}B_3}+ \int_2^3{\cal F}B_3 \\
0  &  \inblue{\int_1^2 (B_4){\cal F}B_2} &   \inblue{\int_1^2B_4{\cal F}B_3} +\inred{\int_2^3(B_4+B_5){\cal F}B_3} \\
0  & 0 &   \inblue{\int_2^3(B_5){\cal F}B_3} \\
0 & 0 & 0 \\
	\end{array} \right. ... \nonumber 
\end{eqnarray}}}
{\tiny{ \begin{eqnarray}
...\left. \begin{array}{llll} 
    \inblue{\int_1^2(B_2){\cal F}B_4} & 0 & 0 & 0\\
\inred{\int_1^2(B_2+B_3){\cal F}B_4}+  \inblue{ \int_2^3(B_3){\cal F}B_4} &  \inblue{\int_2^3(B_3){\cal F}B_5} &   0 & 0\\
    \int_1^2 ){\cal F}B_4  +  \inred{\int_2^3(B_3+B_4){\cal F}B_4} +&   \inblue{\int_3^4(B_4){\cal F}B_6 } + & 0\\
    \quad + \inblue{ \int_3^4(B_4){\cal F}B_4} & + \inblue{\int_3^4(B_4){\cal F}B_6 } & 0 & 0 \\
    \inred{\int_1^2(B_3+B_4){\cal F}B_4} +\int_2^3{\cal F}B_4 +&  \int_2^3{\cal F}B_5+\inred{\int_3^4(B_4+B_5){\cal F}B_5} &  \inred{\int_3^4(B_4+B_5){\cal F}B_6}  & 0  \\
    \quad +\inred{\int_3^4(B_4+B_5){\cal F}B_4} &   &   \\
    \inblue{\int_1^2(B_4){\cal F}B_4}+\inred{\int_2^3(B_4+B_5){\cal F}B_4}+&  \inred{\int_2^3(B_4+B_5){\cal F}B_5}+  \int_3^4{\cal F}B_5+& \inred{\int_2^3(B_4+B_5){\cal F}B_6}+ \int_3^4{\cal F}B_6+  \\
 \quad +\int_1^4{\cal F}B_4 &  \quad + \inred{\int_4^5(B_5+B_6){\cal F}B_5} &\quad + \inred{\int_4^5(B_5+B_6){\cal F}B_6}  &  \\
   \inblue{\int_2^3(B_5){\cal F}B_4} +\inred{\int_3^4(B_5+B_6){\cal F}B_4} &  \inblue{\int_2^3B_5{\cal F}B_5}+ \inred{\int_3^4(B_5+B_6){\cal F}B_5}+&\inblue{\int_2^3B_5{\cal F}B_6}+\inred{\int_3^4(B_5+B_6){\cal F}B_6} + &  \inblue{ \int_2^3B_5{\cal F}B_7} + \\
  &\quad + \int_4^5{\cal F}B_5 & \quad  + \int_4^5{\cal F}B_6 &   + \inred{ \int_3^4(B_5+B_6){\cal F}B_7} + \int_4^5{\cal F}B_7 \\
    \inblue{\int_3^4B_6{\cal F}B_4} + \inred{\int_4^5(B_6+B_7){\cal F}B_4} & \inblue{ \int_3^4B_6{\cal F}B_5 } + \inred{ \int_4^5(B_6+B_7){\cal F}B_5 } &  \inblue{\int_3^4B_6{\cal F}B_6} + \inred{\int_4^5(B_6+B_7){\cal F}B_6}  & \inred{ \int_4^5(B_6+B_7){\cal F}B_7} \\
 	\end{array} \right] \nonumber 
\end{eqnarray}}}

This matrix can be represented as the sum of three sub-matrices

{\footnotesize \begin{eqnarray}
\begin{bmatrix} 
0 & 0 & 0 & 0 & 0 & 0 & 0\\
\int_0^1{\cal F}B_1  &     \int_0^1{\cal F}B_2&     \int_0^1{\cal F}B_3 &  0 & 0 &  0  & 0 \\
 0&     \int_1^2{\cal F}B_2 &     \int_1^2{\cal F}B_3 &    \int_1^2{\cal F}B_4 & 0  &  0 & 0   \\
   0   &  0 &     \int_2^3{\cal F}B_3   & \int_2^3{\cal F}B_4   &     \int_2^3{\cal F}B_5 &  0  & 0\\
   0&     0 & 0 &  \int_3^4{\cal F}B_4  &     \int_3^4{\cal F}B_5  &     \int_3^4{\cal F}B_6 & 0 \\
  0&     0 & 0 &  0 & \int_3^4{\cal F}B_5  &     \int_3^4{\cal F}B_6  &     \int_3^4{\cal F}B_7   \\
0 & 0 & 0 & 0 &    0 &0  & 0 & 0\\
	\end{bmatrix}  + \nonumber
\end{eqnarray} }

{\tiny\begin{eqnarray}
\begin{bmatrix}
\inred{\int_0^1(B_1+B_2){\cal F}B_1}  &     \inred{\int_0^2(B_1+B_2){\cal F}B_2} &  \inred{\int_0^2(B_1+B_2){\cal F}B_3} & 0 & 0 &  0 &  0 \\
0 &  \inred{\int_1^2(B_2+B_3){\cal F}B_2}  &   \inred{\int_1^2(B_3+B_4){\cal F}B_3} &    \inred{\int_1^2(B_3+B_4)B_4} &    0 & 0  & 0\\
     \inred{\int_0^1(B_2+B_3){\cal F}B_1} &    \inred{\int_0^1(B_2+B_3){\cal F}B_2} &  \inred{\int_0^1(B_2+B_3){\cal F}B_3}  +       &     \inred{\int_2^3(B_3+B_4){\cal F}B_4} &  \inred{\int_2^3(B_3+B_4){\cal F}B_5} & 0\\
     & & \quad +     \inred{\int_2^3(B_3+B_4){\cal F}B_3}  &      & & \\
    0  &     \inred{\int_1^2(B_3+B_4){\cal F}B_2} &     \inred{\int_1^2(B_3+B_4){\cal F}B_3} &   \inred{\int_1^2(B_3+B_4){\cal F}B_4} +    &       \inred{\int_3^4(B_4+B_5){\cal F}B_5}  & \inred{\int_3^4(B_4+B_5){\cal F}B_6} & 0\\
      &      &      & \quad +     \inred{\int_3^4(B_4+B_5){\cal F}B_4} &       & \\
   0 & 0  &     \inred{\int_2^3(B_4+B_5){\cal F}B_3} &     \inred{\int_2^3(B_4+B_5){\cal F}B_4} &   \inred{\int_2^3(B_4+B_5){\cal F}B_5} +    &       \inred{\int_4^5(B_5+B_6){\cal F}B_6}  & \inred{\int_4^5(B_5+B_6){\cal F}B_7} \\
&       &      &      & \quad +     \inred{\int_4^5(B_5+B_6){\cal F}B_5} &       & \\
0  &  0 & 0 & \inred{\int_3^4(B_5+B_6){\cal F}B_4}  &   \inred{\int_3^4(B_5+B_6){\cal F}B_5} &     \inred{\int_3^4(B_5+B_6){\cal F}B_6}  &  0  \\
    0  &   0 &  0 & 0 &   \inred{\int_4^5(B_6+B_7){\cal F}B_5} & \inred{\int_4^5(B_6+B_7){\cal F}B_6} &     \inred{\int_4^5(B_6+B_7){\cal F}B_7}  \\
	\end{bmatrix}+ \nonumber
\end{eqnarray} }

{\footnotesize \begin{eqnarray}
\begin{bmatrix} 
0 & 0 &     0 & 0 & 0 &    0 & 0\\
0 &  0 &  \inblue{\int_2^3B_3{\cal F}B_3}  &    \inblue{\int_2^3B_3{\cal F}B_4} & \inblue{\int_2^3B_3{\cal F}B_5}   &  0 & 0\\
    0 &   0 & 0 &  \inblue{\int_3^4B_4{\cal F}B_4} &  \inblue{\int_3^4B_4{\cal F}B_5}   &  \inblue{\int_3^4B_4{\cal F}B_6} & 0\\
  \inblue{\int_0^1B_3{\cal F}B_1}  &      \inblue{\int_0^1B_3{\cal F}B_2} &  \inblue{\int_0^1B_3{\cal F}B_3}  & 0 &   \inblue{\int_4^5B_5{\cal F}B_5} &  \inblue{\int_4^5B_5{\cal F}B_6}   &  \inblue{\int_4^5B_5{\cal F}B_7} \\
   0  &     \inblue{\int_1^2B_4{\cal F}B_2} &     \inblue{\int_1^2B_4{\cal F}B_3}   &   \inblue{\int_1^2B_4{\cal F}B_4} & 0 & 0 & 0\\
  0&   0  &     \inblue{\int_2^3B_5{\cal F}B_3} &     \inblue{\int_2^3B_5{\cal F}B_4} & \inblue{\int_2^3B_5{\cal F}B_5} &  0   & 0 \\
 0 &  0&   0  &     \inblue{\int_3^4B_6{\cal F}B_4} &     \inblue{\int_3^4B_6{\cal F}B_5} & \inblue{\int_3^4B_6{\cal F}B_6} & 0    \\
	\end{bmatrix} \nonumber
\end{eqnarray} }

The black terms, represents the case, where we have the summation of all local test functions over a single element.  
In such the case

$\sum_{m=i,i+p}\int_E{B^x_{m,p}{\cal F}\left(B^x_{j,p}\right)} $ =
$ \sum_o w_o{[\sum_{m=i,i+p} B^x_{m,p}(x_o)] {\cal F}\left(B^x_{j,p}(x_o)\right)}Jac(x_o)$ = 

$\sum_o w_o{{\cal F}\left(B^x_{j,p}(x_o)\right)}Jac(x_o) = 1\int_E{{\cal F}\left(B^x_{j,p}\right)} $

since the test function sum up to one $[\sum_{m=j,j+p} B^x_{m,p}(x_o)]=1$, from the partition of unity property.
The black terms represent the test functions equal to 1.

The red terms represent the case, where we have the integration over a single element of a sum of two test B-splines multiplied by our operator applied to a trial function.

The blue terms represent the integration over a single element with one test B-spline multiplied by our operator applied to a trial function.

The blue and red terms, they cannot be removed from the system. However, we will show how to make their contribution negligible.
Their presence in a matrix is a consequence of the fact that we sum up three B-splines that span over different three elements, and over the beginning and the last two elements, they do not sum up to one. They only sum up to one over the central element.

\begin{figure}
\includegraphics[scale=0.5]{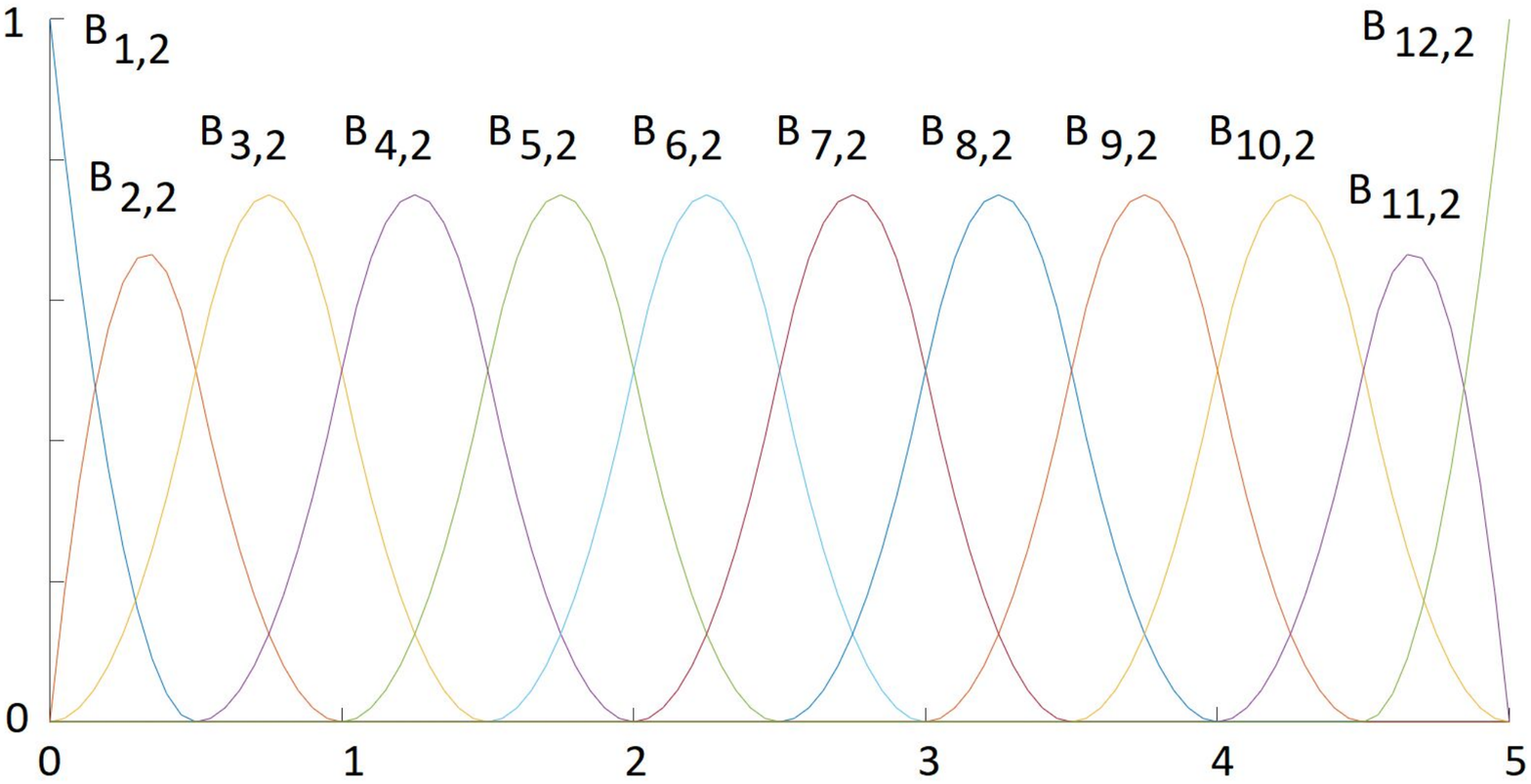}
\caption{B-splines span over [0 0 0 0.5 1 1.5 2 2.5 3 3.5 4 4.5 5 5 5] knot vector.}
\label{fig::BsplinesC1adapted}
\end{figure}

Let us check what happens if we sum up more test B-splines, and increase the number of elements over the test space only.
Let us double the number of elements over the test space, by taking the knot vector [0 0 0 0.5 1 1.5 2 2.5 3 3.5 4 4.5 5 5 5]. We have now the test functions presented in Figure \ref{fig::BsplinesC1adapted}.
If we sum up three test B-splines, we will get a single central segment where the test B-splines sum up to one, this time thinner, since the refined elements are smaller than the original elements.
If we sum up more rows of the system, we will get a longer interval where B-splines sum up to one.
For the sum of four rows of the matrix, representing four test B-splines we get the function constant on the central segment [2 3], and the "blue" and "red" terms they become two times smaller since the corresponding segments of B-splines are two times "thinner".

\begin{figure}
\includegraphics[scale=0.5]{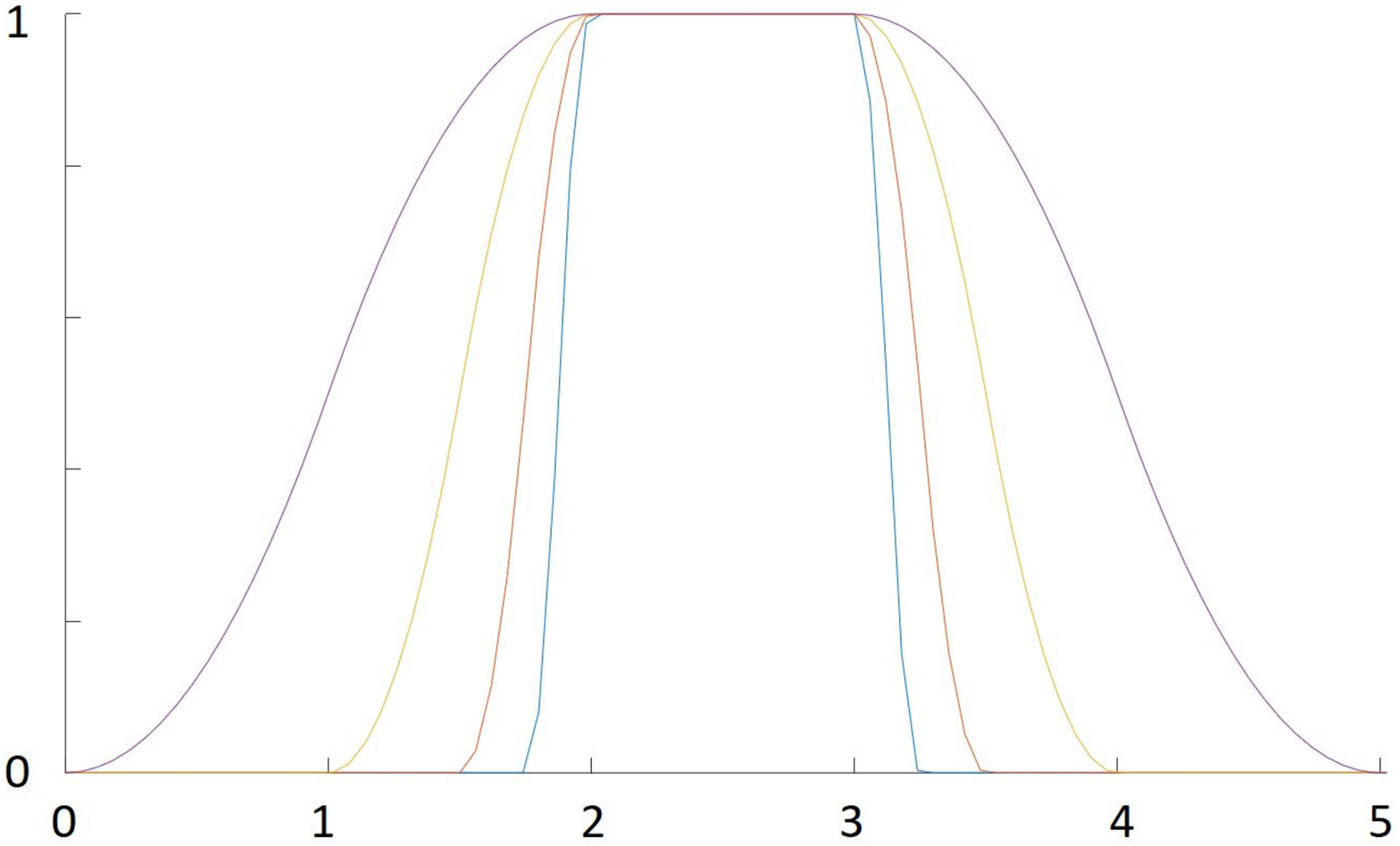}
\caption{The convergence of the test function to the piece-wise constant functions, when we refine the mesh and increase the number of summed up rows.}
\label{fig:convergence}
\end{figure}

Increasing the number of test B-splines further, and summing up more test B-splines results in convergence to the piece-wise constant test functions, as presented in Figure \ref{fig:convergence}. By changing the range of the summation of test B-splines, we can change the location of the segment. We can cover any interval of elements by a segment where the resulting test function is equal to 1. We have extra two thin segments at the beginning and at the end of the constant segment, where the shape is smoothly increasing from 0 to 1.

For a given mesh, we can sum up sets of three B-splines and get test functions with one segment equal to one, and the two other segments being quadratic polynomials.
We can also sum sets of more functions and get "longer" segments equal to one, again with the two segments, at the beginning and at at the end, being quadratic polynomials.

Now, the question is, how to get rid of the polynomial segments at the beginning and at the end of the test functions? How to work with piece-wise constant test functions?
When we increase the number of elements and the length of the segments equal to 1, the contribution of ``red'' and ``blue'' terms become negligibly small. At the limit (when we increase the number of elements and number of added test functions), they vanish.

We construct the isogeometric analysis method with piece-wise constant test functions in the following way
\begin{itemize}
\item We fix the trial space, with the trial B-spline basis functions $\{ B^x_{i,p}(x)\}_{i=1,...,N_x}$.
\item We plug our trial B-splines into our PDE, namely $u={\cal RHS}$, so we have $\sum_{i=1,...,N_x} u_{i,j}B^x_{i,p}(x) = RHS(x)$.
\item We take the test space $\{\hat{B}^x_{j,p}\}_{j=1,...,N^{*}_x}$, larger than trial space, with $N^{*}_x>>N_x$, and we multiply our equation and integrate. In other words, we take scalar L2 products with more test functions than trial functions.
\begin{eqnarray}
\int  \hat{B}^x_{j,p} \sum_{i=1,...,N_x} u_{i}{\cal F } B^x_{i,p}(x) dx = \int \hat{B}^x_{j,p}  RHS(x) dx \quad j=1,...,N^{*}_x>>N_x \nonumber
\end{eqnarray}
We end up with the rectangular matrix
\begin{eqnarray}
\begin{bmatrix}
\int  \hat{B}^x_{1,p} {\cal F } B^x_{1,p}(x) dx  & \cdots & \int  \hat{B}^x_{1,p} {\cal F } B^x_{N_x,p}(x) dx  \\
\int  \hat{B}^x_{2,p} {\cal F } B^x_{1,p}(x) dx  & \cdots & \int  \hat{B}^x_{2,p} {\cal F } B^x_{N_x,p}(x) dx  \\
\vdots & \vdots & \vdots \\
\int  \hat{B}^x_{N^{*}_x,p} {\cal F } B^x_{1,p}(x) dx  & \cdots & \int  \hat{B}^x_{N^{*}_x,p} {\cal F } B^x_{N_x,p}(x) dx  \\
\end{bmatrix}
\begin{bmatrix}
u_1 \\
u_2 \\
\vdots \\
u_{N_x} \\
\end{bmatrix}
=  \begin{bmatrix}
\int \hat{B}^x_{1,p}  RHS(x) dx  \\
\int \hat{B}^x_{2,p}  RHS(x) dx \\
\vdots \\
\int \hat{B}^x_{N^{*}_x,p}  RHS(x) dx  \\
\end{bmatrix} \nonumber
\end{eqnarray}
\item Now, select $N_x$ sets of indices $\{{\cal J}_i\}_{i=1,...,N_x}, {\cal J}_i \subset \{1,...,N^{*}_x\}$, and we sum up the corresponding equations into the new system.
\begin{eqnarray}
\begin{bmatrix}
\int \sum_{m\in {\cal J}_1} \hat{B}^x_{m,p} {\cal F } B^x_{1,p}(x) dx  & \cdots & \int  \sum_{m\in {\cal J}_1} \hat{B}^x_{m,p} {\cal F } B^x_{N_x,p}(x) dx  \\
\int \sum_{m\in {\cal J}_2} \hat{B}^x_{m,p} {\cal F } B^x_{1,p}(x) dx  & \cdots & \int  \sum_{m\in {\cal J}_2}\hat{B}^x_{m,p} {\cal F } B^x_{N_x,p}(x) dx  \\
\vdots & \vdots & \vdots \\
\int  \sum_{m\in {\cal J}_{N_x}}\hat{B}^x_{m,p} {\cal F } B^x_{1,p}(x) dx  & \cdots & \int  \sum_{m\in {\cal J}_{N_x}} \hat{B}^x_{m,p} {\cal F } B^x_{N_x,p}(x) dx  \\
\end{bmatrix}
\begin{bmatrix}
u_1 \\
u_2 \\
\vdots \\
u_{N_x} \\
\end{bmatrix}
= \nonumber \\
\begin{bmatrix}
\int  \sum_{m\in {\cal J}_1} \hat{B}^x_{m,p}  RHS(x) dx  \\
\int  \sum_{m\in {\cal J}_2} \hat{B}^x_{m,p}  RHS(x)  \\
\vdots \\
\int \sum_{m\in {\cal J}_{N_x}} \hat{B}^x_{m,p}  RHS(x) dx  \\
\end{bmatrix} \nonumber
\end{eqnarray}

We do it in such a way that the obtained system is well-posed (that the linear combinations of test functions from the selected subsets of test functions form a linearly independent basis).
Namely, we select the intervals over our domain, where we want our piece-wise constant test functions to be fixed to one. We select and sum up rows in such a way, that we end up with piece-wise constant test functions $\{{\cal I}_i\}_{i=1,...,N_x}$ span over some intervals. We select intervals in such a way that they are not linearly dependent to the obtained well-posed system of equations.
\end{itemize}

We end up with the system of equations

\begin{eqnarray}
	\begin{bmatrix}
    \int {\cal I}_{1}{{\cal F}\left(B^x_{1,p}\right)}dx  & \cdots &  \int {\cal I}_{1}{{\cal F}\left(B^x_{N_x,p}\right)}dx \\
    \int {\cal I}_{2}{{\cal F}\left(B^x_{1,p}\right)}dx & \cdots &  \int {\cal I}_{2} {{\cal F}\left(B^x_{N_x,p}\right)}dx\\   
    \vdots & \vdots &  \vdots \\
\int {\cal I}_{N_x}{{\cal F}\left(B^x_{1,p}\right)}dx  & \cdots &  \int {\cal I}_{N_x}{{\cal F}\left(B^x_{N_x,p}\right)}dx \\
  \end{bmatrix} 
\begin{bmatrix}
  u_1 \\ u_2\\  \vdots \\ u_{N_{x}}  \\
  \end{bmatrix} 
=    \begin{bmatrix}
\int {\cal I}_{1} {\cal RHS}dx \\
\int {\cal I}_{2}{\cal RHS}dx \\
\vdots \\
\int {\cal I}_{N_x} {\cal RHS} dx\\
  \end{bmatrix}  \nonumber
\end{eqnarray}

The considerations for higher-order B-splines follows similar lines as for the quadratic B-splines.

In general, summing $2p+1$ B-splines of order $p$, gives the test function over one element equal to 1, and over $p$ elements at the beginning, and $p$ elements at the end, where the test functions change smoothly from 0 to 1.
Summing $2p+m$ B-splines of order $p$, gives test functions equal to 1 over $m$ elements, and $p$ segments at the beginning and at the end, where the function is smoothly going from 0 to 1.
In the limit, using more elements over the test space, and summing more rows, we can get a segment equal to 1 over any interval span over our trial space.

Selecting the piece-wise constant test functions has to be done in such a way that they are linearly independent, and the resulting system of equations can be factorized using direct solver. We must select intervals in such a way that the number of test functions is equal to the number of trial functions, and the test functions are linearly independent. Otherwise, the factorization will break.

\section{Two dimensional case}

We repeat our considerations in the two-dimensional case. 
We start from the general form of the PDE
\begin{equation}
{\cal F}u={\cal RHS}
\end{equation}
where we discretize with $C^1$ continuity B-splines, and we do not integrate by parts.
We have the global system of linear equations
{\footnotesize\begin{equation*}
	\begin{bmatrix}
    \int{\left(B^x_{1,p}B^y_{1,p}\right)}{{\cal F}\left(B^x_{1,p}B^y_{1,p}\right)} & \int{\left(B^x_{1,p}B^y_{1,p}\right)}{{\cal F}\left(B^x_{2,p}B^y_{1,p}\right)} & \cdots &  \int{\left(B^x_{1,p}B^y_{1,p}\right)}{{\cal F}\left(B^x_{N_x,p}B^y_{N_y,p}\right)} \\
    \int{\left(B^x_{2,p}B^y_{1,p}\right)}{{\cal F}\left(B^x_{1,p}B^y_{1,p}\right)} & \int{\left(B^x_{2,p}B^y_{1,p}\right)}{{\cal F}\left(B^x_{2,p}B^y_{1,p}\right)} & \cdots &  \int{\left(B^x_{2,p}B^y_{1,p}\right)}{{\cal F}\left(B^x_{N_x,p}B^y_{N_y,p}\right)} \\   
    \vdots & \vdots & \vdots &  \vdots \\
        \int{\left(B^x_{N_x,p}B^y_{N_y,p}\right)}{{\cal F}\left(B^x_{1,p}B^y_{1,p}\right)} & \int{\left(B^x_{N_x,p}B^y_{N_y,p}\right)}{{\cal F}\left(B^x_{2,p}B^y_{1,p}\right)} & \cdots &  \int{\left(B^x_{N_x,p}B^y_{N_y,p}\right)}{{\cal F}\left(B^x_{N_x,p}B^y_{N_y,p}\right)} \\
  \end{bmatrix}
  \begin{bmatrix}
  {u_{1,1}} \\ {u_{2,1}} \\ \vdots \\ {u_{N_x,N_y}}\\
  \end{bmatrix} 
\end{equation*}
\begin{equation*} 
=    \begin{bmatrix}
\int {\cal RHS}(x,y) {B^x_{1,p}B^y_{1,p}}  \\
\int {\cal RHS}(x,y) {B^x_{2,p}B^y_{1,p}}  \\
\vdots \\
\int {\cal RHS}(x,y) {B^x_{N_x,p}B^y_{N_y,p}} \\
  \end{bmatrix}
\end{equation*}}

We consider a quadrature with points and weights $\{(x_o,y_o),w_o\}_o$.
At a given point $(x_o,y_o)$ from the selected quadrature, we have $2p-1$ non-zero B-spline basis functions in one direction.

Since each row in the global matrix corresponds to one test function $B^x_{i,p}B^y_{j,p}B^x_{m,p}B^y_{n,p}$ we can number them $(i,j;m,n)$.

We select the $N^{*}_e$ intervals of the test functions along $x$ direction, we adapt the test space in the $x$ direction, and sum up with multiple rows of the test space, to get the piece-wise constant test functions over the selected intervals.

Namely, we add to the row $(i,j;m,n)$ the sum of rows number 

\begin{equation}
(min(1,i-k),j;m,n),...,(max(N^{*}_x,i+k),j;m,n)
\end{equation}

where $N^{*}_x=N^{*}_e+p$ denotes the number of test functions in the $x$ direction.

We get the equivalent global system 

{\footnotesize\begin{equation*}
	\begin{bmatrix}
   \sum_{m=1,...,k+1}   \int{\left(B^x_{m,p}B^y_{1,p}\right)}{{\cal F}\left(B^x_{1,p}B^y_{1,p}\right)} & \cdots &    \sum_{m=1,...,k+1} \int{\left(B^x_{m,p}B^y_{1,p}\right)}{{\cal F}\left(B^x_{N_x,p}B^y_{N_y,p}\right)} \\
\vdots & \vdots & \vdots  \\
     \sum_{m=i-k,...,i+k}  \int{\left(B^x_{m,p}B^y_{j,p}\right)}{{\cal F}\left(B^x_{1,p}B^y_{1,p}\right)} & \cdots &        \sum_{m=i-k,...,i+k} \int{\left(B^x_{m,p}B^y_{j,p}\right)}{{\cal F}\left(B^x_{N_x,p}B^y_{N_y,p}\right)} \\   
    \vdots & \vdots & \vdots \\
\sum_{m=N^{*}_x-k,...,N^{*}_x}        \int{\left(B^x_{m,p}B^y_{N_y,p}\right)}{{\cal F}\left(B^x_{1,p}B^y_{1,p}\right)} & \cdots &  \sum_{m=N^{*}_x-k,...,N^{*}_x}    \int{\left(B^x_{m,p}B^y_{N_y,p}\right)}{{\cal F}\left(B^x_{N_x,p}B^y_{N_y,p}\right)} \\
  \end{bmatrix}
  \begin{bmatrix}
  {u_{1,1}} \\ \vdots\\ {u_{i,j}} \\ \vdots \\ {u_{N_x,N_y}}\\
  \end{bmatrix} 
\end{equation*}
\begin{equation*} 
=    \begin{bmatrix}
\sum_{m=1,...,k+1}\int {\cal RHS}(x,y) {B^x_{m,p}(x)B^y_{1,p}(y)}  \\
\vdots \\
\sum_{m=i-k,...,i+k}\int {\cal RHS}(x,y) {B^x_{m,p}(x)B^y_{j,p}(y)}  \\
\vdots \\
\sum_{m=N^{*}_x-k,...,N^{*}_x}\int {\cal RHS}(x,y) {B^x_{m,p}(x)B^y_{N_y,p}(y)} \\
  \end{bmatrix}
\end{equation*}}

Now, we compute the integrals by using numerical quadrature rule for polynomials

\begin{equation*}
\begin{bmatrix}
  \sum_o w_o{\left([\sum_{m=1,...,k+1}  B^x_{m,p}(x_o)] B^y_{1,p}(y_o)\right)}{{\cal F}\left(B^x_{1,p}(x_o)B^y_{1,p}(y_o)\right)}Jac(x_o,y_o)  & \cdots 
\\
    \vdots & \vdots &  \\
  \sum_o w_o{\left([\sum_{m=i-k,...,i+k}  B^x_{m,p}(x_o)] B^y_{j,p}(y_o)\right)}{{\cal F}\left(B^x_{1,p}(x_o)B^y_{1,p}(y_o)\right)}Jac(x_o,y_o)  & \cdots 
  \\  \vdots & \vdots \\
        \sum_o w_o{\left([\sum_{m=N^{*}_x-k,...,N^{*}_x}  B^x_{m,p}(x_o)]B^y_{N_y,p}(y_o)\right)}{{\cal F}\left(B^x_{1,p}(x_o)B^y_{1,p}(y_o)\right)}Jac(x_o,y_o) & \cdots 
\\
  \end{bmatrix}
  \begin{bmatrix}
  {u_{1,1}} \\ \vdots \\ {u_{i,j}}\\ \vdots \\ {u_{N_x,N_y}}\\
  \end{bmatrix} 
\end{equation*}
\begin{equation*} 
=    \begin{bmatrix}
\sum_o w_o{\cal RHS}(x_o,y_o) {[\sum_{m=1,...,k+1}B^x_{m,p}(x_o)]B^y_{1,p}(y_o)}Jac(x_o,y_o) \\
\vdots \\
\sum_o w_o{\cal RHS}(x_o,y_o) {[\sum_{m=i-k,...,i+k}B^x_{m,p}(x_o)]B^y_{j,p}(y_o)}Jac(x_o,y_o) \\
\vdots \\
\sum_o w_o{\cal RHS}(x_o,y_o) {[\sum_{m=N^{*}_x-k,...,N^{*}_x}B^x_{m,p}(x_o)]B^y_{N_y,p}(y_o)}Jac(x_o,y_o)  \\
  \end{bmatrix}
\end{equation*}

At a given quadrature point, we sum up all the B-splines in one direction.
The number of test functions that we sum up at a given row is such that the summation, from the partition of unity, is equivalent to the piece-wise constant test function in the $x$ direction.
The other terms (the ``blue'' and the ``red'' terms) they are neglected (or they disappear in the limit).

So these summation terms disappear.
\begin{equation*}
\begin{bmatrix}
  \sum_o w_o{\left( B^y_{1,p}(y_o)\right)}{{\cal F}\left(B^x_{1,p}(x_o)B^y_{1,p}(y_o)\right)}Jac(x_o,y_o)  & \cdots
\\
    \vdots & \vdots &  \\
  \sum_o w_o{\left(B^y_{j,p}(y_o)\right)}{{\cal F}\left(B^x_{1,p}(x_o)B^y_{1,p}(y_o)\right)}Jac(x_o,y_o)  & \cdots 
  \\  \vdots & \vdots \\
        \sum_o w_o{\left(B^y_{N_y,p}(y_o)\right)}{{\cal F}\left(B^x_{1,p}(x_o)B^y_{1,p}(y_o)\right)}Jac(x_o,y_o) & \cdots
\\
  \end{bmatrix}
  \begin{bmatrix}
  {u_{1,1}} \\ \vdots \\ {u_{i,j}}\\ \vdots \\ {u_{N_x,N_y}}\\
  \end{bmatrix} 
\end{equation*}
\begin{equation*} 
=    \begin{bmatrix}
\sum_o w_o{\cal RHS}(x_o,y_o) {B^y_{1,p}(y_o)}Jac(x_o,y_o) \\
\vdots \\
\sum_o w_o{\cal RHS}(x_o,y_o) {B^y_{j,p}(y_o)}Jac(x_o,y_o) \\
\vdots \\
\sum_o w_o{\cal RHS}(x_o,y_o) {B^y_{N_y,p}(y_o)}Jac(x_o,y_o)  \\
  \end{bmatrix}
\end{equation*}

Now, we can come back to the integral, and we have now the piece-wise constant test functions $I_{i}(x)$.

\begin{equation*}
	\begin{bmatrix}
  \int {\cal I}_{1}B^y_{1,p}{{\cal F}\left(B^x_{1,p}B^y_{1,p}\right)}  & \cdots &  \int {\cal I}_{1}B^y_{1,p}{{\cal F}\left(B^x_{N_x,p}B^y_{N_y,p}\right)} \\  
    \vdots & \vdots &  \vdots \\
  \int {\cal I}_{i}B^y_{j,p}{{\cal F}\left(B^x_{1,p}B^y_{1,p}\right)}  & \cdots &  \int {\cal I}_{i}B^y_{j,p}{{\cal F}\left(B^x_{N_x,p}B^y_{N_y,p}\right)} \\  
    \vdots & \vdots &  \vdots \\
        \int {\cal I}_{N_x}B^y_{N_y,p}{{\cal F}\left(B^x_{1,p}B^y_{1,p}\right)} & \cdots &  \int {\cal I}_{N_x}B^y_{N_y,p}{{\cal F}\left(B^x_{N_x,p}B^y_{N_y,p}\right)} \\
  \end{bmatrix}
  \begin{bmatrix}
  {u_{1,1}} \\ \vdots \\ {u_{i,j}}\\ \vdots \\ {u_{N_x,N_y}}\\
  \end{bmatrix} 
\end{equation*}

\begin{equation*} 
=    \begin{bmatrix}
\int {\cal I}_{1}{\cal RHS}(x,y) {B^y_{1,p}(y)} \\
\vdots \\
\int {\cal I}_{i}{\cal RHS}(x,y) {B^y_{j,p}(y)} \\
\vdots \\
\int {\cal I}_{N_x}{\cal RHS}(x,y) {B^y_{N_y,p}(y)}  \\
  \end{bmatrix}
\end{equation*}

Now, we repeat the same logic with respect to the one-dimensional B-spline basis functions in the $y$ direction.

We select the $N^{*}_y$ elements of the test functions along $y$ direction, we adapt the test space in the $y$ direction, and sum up with multiple rows of the test space, to get the piece-wise constant test functions over the selected intervals.
We get

\begin{equation*}
	\begin{bmatrix}
  \int {\cal I}_{1}{\cal I}_{1}{{\cal F}\left(B^x_{1,p}B^y_{1,p}\right)}  & \cdots &  \int {\cal I}_{1}{\cal I}_{1}{{\cal F}\left(B^x_{N_x,p}B^y_{N_y,p}\right)} \\  
    \vdots & \vdots &  \vdots \\
  \int {\cal I}_{i}{\cal I}_{j}{{\cal F}\left(B^x_{1,p}B^y_{1,p}\right)}  & \cdots &  \int {\cal I}_{i}{\cal I}_{j}{{\cal F}\left(B^x_{N_x,p}B^y_{N_y,p}\right)} \\  
    \vdots & \vdots &  \vdots \\
        \int {\cal I}_{N_x}{\cal I}_{N_y}{{\cal F}\left(B^x_{1,p}B^y_{1,p}\right)} & \cdots &  \int {\cal I}_{N_x}{\cal I}_{N_y}{{\cal F}\left(B^x_{N_x,p}B^y_{N_y,p}\right)} \\
  \end{bmatrix}
  \begin{bmatrix}
  {u_{1,1}} \\ \vdots \\ {u_{i,j}}\\ \vdots \\ {u_{N_x,N_y}}\\
  \end{bmatrix} 
\end{equation*}

\begin{equation*} 
=    \begin{bmatrix}
\int {\cal I}_{1}{\cal I}_{1}{\cal RHS}(x,y) \\
\vdots \\
\int {\cal I}_{i}{\cal I}_{j}{\cal RHS}(x,y) \\
\vdots \\
\int {\cal I}_{N_x}{\cal I}_{N_y}{\cal RHS}(x,y)  \\
  \end{bmatrix}
\end{equation*}

\section{Examples}

In this section, we present four numerical examples. The goal of the first example is to show how the method scales on a three-dimensional projection problem if we increase the mesh size or the B-splines order. The goal of the second example is to illustrate that the method can be applied for explicit dynamics problems since each of them is a sequence of isogeometric L2 projections. The goal of the third example is to show that the method allows incorporating boundary conditions. We also show how the method scales with a two-dimensional MATLAB code. Finally, we show the comparison of our method with the isogeometric L2 projection of a bitmap. We compare the convergence rates on this difficult projection example.   

\subsection{Isogeometric L2 projections}
First example is the L2 projection problem.
\begin{equation}
u={\cal RHS} \nonumber
\end{equation}
which in the weak form is
\begin{equation}
(v,u)=(v,{\cal RHS}) \nonumber
\end{equation}
solved over $\Omega=[0,1]^3$.

We define the B-spline basis for trial and test $\{ B^x_{i,p}B^y_{j,p}B^z_{k,p} \}_{i=1,...,N_x;j=1,...,N_y;k=1,...,N_z}$ and we discretize in the standard Galerkin way

{\footnotesize \begin{equation*}
	\begin{bmatrix}
    \int{\left(B^x_{1,p}B^y_{1,p}B^z_{1,p}\right)}{\left(B^x_{1,p}B^y_{1,p}B^z_{1,p}\right)} & \cdots &  \int{\left(B^x_{1,p}B^y_{1,p}B^z_{1,p}\right)}{\left(B^x_{N_x,p}B^y_{N_y,p}B^z_{N_z,p}\right)} \\
    \int{\left(B^x_{2,p}B^y_{1,p}B^z_{1,p}\right)}{\left(B^x_{1,p}B^y_{1,p}B^z_{1,p}\right)} & \cdots &  \int{\left(B^x_{2,p}B^y_{1,p}B^z_{1,p}\right)}{\left(B^x_{N_x,p}B^y_{N_y,p}B^z_{N_z,p}\right)} \\   
    \vdots & \vdots & \vdots \\
        \int{\left(B^x_{N_x,p}B^y_{N_y,p}B^z_{N_z,p}\right)}{\left(B^x_{1,p}B^y_{1,p}B^z_{1,p}\right)}& \cdots &  \int{\left(B^x_{N_x,p}B^y_{N_y,p}B^z_{N_z,p}\right)}{\left(B^x_{N_x,p}B^y_{N_y,p}B^z_{N_z,p}\right)} \\
  \end{bmatrix}
  \begin{bmatrix}
  {u_{1,1,1}} \\ {u_{2,1,1}} \\ \vdots \\ {u_{N_x,N_y,N_z}}\\
  \end{bmatrix} 
\end{equation*}
\begin{equation*} 
=    \begin{bmatrix}
\int {\cal RHS}(x,y,z) {B^x_{1,p}(x)B^y_{1,p}(y)B^z_{1,p}(z)}  \\
\int {\cal RHS}(x,y,z) {B^x_{2,p}(x)B^y_{1,p}(y)B^z_{1,p}(z)}  \\
\vdots \\
\int {\cal RHS}(x,y,z) {B^x_{N_x,p}(x)B^y_{N_y,p}(y)B^z_{N_z,p}(z)} \\
  \end{bmatrix}
\end{equation*}}

Now, we can set the test functions to piece-wise constants and adjust the integrals accordingly to the spans of the test functions

{\footnotesize \begin{equation*}
	\begin{bmatrix}
    \int {\cal I}_{1}{\cal I}_{1}{\cal I}_{1}{\left(B^x_{1,p}B^y_{1,p}B^z_{1,p}\right)} & \cdots &  \int {\cal I}_{1}{\cal I}_{1}{\cal I}_{1}{\left(B^x_{N_x,p}B^y_{N_y,p}B^z_{N_z,p}\right)} \\
    \int {\cal I}_{2}{\cal I}_{1}{\cal I}_{1}{\left(B^x_{1,p}B^y_{1,p}B^z_{1,p}\right)}& \cdots &  \int {\cal I}_{2}{\cal I}_{1}{\cal I}_{1}{\left(B^x_{N_x,p}B^y_{N_y,p}B^z_{N_z,p}\right)} \\   
    \vdots & \vdots & \vdots &  \vdots \\
        \int {\cal I}_{N_x}{\cal I}_{N_y}{\cal I}_{N_z}{\left(B^x_{1,p}B^y_{1,p}B^z_{1,p}\right)} & \cdots &  \int {\cal I}_{N_x}{\cal I}_{N_y}{\cal I}_{N_z}{\left(B^x_{N_x,p}B^y_{N_y,p}B^z_{N_z,p}\right)} \\
  \end{bmatrix}
  \begin{bmatrix}
  {u_{1,1,1}} \\ {u_{2,1,1}} \\ \vdots \\ {u_{N_x,N_y,N_z}}\\
  \end{bmatrix} 
\end{equation*}
\begin{equation*} 
=    \begin{bmatrix}
\int {\cal I}_{1}{\cal I}_{1}{\cal I}_{1} {\cal RHS}(x,y,z)  dxdydz \\
\int {\cal I}_{2}{\cal I}_{1}{\cal I}_{1} {\cal RHS}(x,y,z)  dxdydz \\
\vdots \\
\int {\cal I}_{N_x}{\cal I}_{N_y}{\cal I}_{N_z} {\cal RHS}(x,y,z)  dxdydz \\
  \end{bmatrix}
\end{equation*}}

Let us test the scalability of our method, using standard Gaussian quadrature. We assume that the right-hand side is the polynomial of the third order with respect to each variable, e.g.,
\begin{equation}
{\cal RHS}(x,y,z)=(a_xx^3+b_xx^2+c_xx+d_x)(a_yy^3+b_yy^2+c_yy+d_y)(a_zz^3+b_zz^2+c_zz+d_z) \nonumber
\end{equation}
Standard isogeometric L2 projection for second order B-splines with $C^1$ continuity
\begin{equation}
\int{\cal RHS}(x,y,z)B^x_{i,2}B^y_{j,2}B^z_{k,2}dxdydz = {\cal O}(x^5){\cal O}(y^5){\cal O}(z^5) \nonumber
\end{equation}
requires the third order quadrature, to integrate the right-hand side exactly, since $2*3-1=5$ and we have to integrate polynomials of the fifth order in each direction.

When we introduce piece-wise constant test polynomials,
\begin{equation}
\int{\cal RHS}(x,y,z)dxdydz = {\cal O}(x^3){\cal O}(y^3){\cal O}(z^3) \nonumber
\end{equation}
the exact right-hand side integration requires the second order quadrature, since $2*2-1=3$ and we have to integrate polynomials of the third order in each direction.

We use alternating directions direct solver for factorization \cite{CPC}. This implementation of the direct solver for isogeometric L2 projections has the following features. It has a linear computational cost ${\cal O}(N)$, and it uses the Kronecker product structure of the matrix. It generates three one-dimensional systems with multiple RHS. In the case of piece-wise constant test functions, these systems look in the following way. The first system
\begin{equation*}
	\begin{bmatrix}
    \int {\cal I}_{1}{ B^x_{1,p}}dx & \cdots &  \int {\cal I}_{1}{B^x_{N_x,p}}dx \\
    \int {\cal I}_{2}{B^x_{1,p}}dx & \cdots &  \int {\cal I}_{2}{B^x_{N_x,p}}dx\\   
    \vdots & \vdots & \vdots \\
\int {\cal I}_{N_x}{B^x_{1,p}}dx  
 & \cdots &  \int {\cal I}_{N_x}{B^x_{N_x,p}}dx \\
  \end{bmatrix} 
  \begin{bmatrix}
  {v_{1,1,1}} & \vdots & {v_{1,N_y,N_z}} \\ {v_{2,1,1}} & \vdots & {v_{2,N_y,N_z}} \\ \vdots \\ {v_{N_x,1,1}} & \vdots & {v_{N_x,N_y,N_z}}\\
  \end{bmatrix} 
\end{equation*}
\begin{equation*} 
=    \begin{bmatrix}
\int {\cal I}_{1}{\cal I}_{1}{\cal I}_{1} {\cal RHS}(x,y,z)  dxdydz & \vdots & 
\int {\cal I}_{1}{\cal I}_{N_y}{\cal I}_{N_z} {\cal RHS}(x,y,z)  dxdydz \\
\int {\cal I}_{2}{\cal I}_{1}{\cal I}_{1} {\cal RHS}(x,y,z)  dxdydz & \vdots 
& \int {\cal I}_{2}{\cal I}_{N_y}{\cal I}_{N_z} {\cal RHS}(x,y,z)  dxdydz \\
\vdots \\
\int {\cal I}_{N_x}{\cal I}_{1}{\cal I}_{1} {\cal RHS}(x,y,z)  dxdydz  & \vdots & 
\int {\cal I}_{N_x}{\cal I}_{N_y}{\cal I}_{N_z} {\cal RHS}(x,y,z)  dxdydz  \\
  \end{bmatrix}
\end{equation*}

the second system
\begin{equation*}
	\begin{bmatrix}
    \int {\cal I}_{1}{ B^y_{1,p}}dy & \cdots &  \int {\cal I}_{1}{B^y_{N_x,p}}dy \\
    \int {\cal I}_{2}{B^y_{1,p}}dy & \cdots &  \int {\cal I}_{2}{B^y_{N_x,p}}dy\\   
    \vdots & \vdots & \vdots \\
\int {\cal I}_{N_y}{B^y_{1,p}}dy  
 & \cdots &  \int {\cal I}_{N_z}{B^y_{N_z,p}}dy \\
  \end{bmatrix} 
  \begin{bmatrix}
  {w_{1,1,1}} & \vdots & {w_{N_x,1,N_z}} \\ {w_{1,2,1}} & \vdots & {w_{N_x,2,N_z}} \\ \vdots \\ {w_{1,N_y,1}} & \vdots & {w_{N_x,N_y,N_z}}\\
  \end{bmatrix} 
=
  \begin{bmatrix}
  {v_{1,1,1}} & \vdots & {v_{N_x,1,N_z}} \\ {v_{1,2,1}} & \vdots & {v_{N_x,2,N_z}} \\ \vdots \\ {v_{1,N_y,1}} & \vdots & {v_{N_x,N_y,N_z}}\\
  \end{bmatrix} 
\end{equation*}

and the third system
\begin{equation*}
	\begin{bmatrix}
    \int {\cal I}_{1}{ B^z_{1,p}}dz & \cdots &  \int {\cal I}_{1}{B^z_{N_z,p}}dz \\
    \int {\cal I}_{2}{B^z_{1,p}}dz & \cdots &  \int {\cal I}_{2}{B^z_{N_z,p}}dz\\   
    \vdots & \vdots & \vdots \\
\int {\cal I}_{N_z}{B^z_{1,p}}dz  
 & \cdots &  \int {\cal I}_{N_z}{B^z_{N_z,p}}dz \\
  \end{bmatrix} 
  \begin{bmatrix}
  {u_{1,1,1}} & \vdots & {u_{N_x,N_y,1}} \\ {u_{1,1,2}} & \vdots & {u_{N_x,N_y,2}} \\ \vdots \\ {u_{1,1,N_z}} & \vdots & {w_{N_x,N_y,N_z}}\\
  \end{bmatrix} 
=
  \begin{bmatrix}
  {w_{1,1,1}} & \vdots & {w_{N_x,N_y,1}} \\ {w_{1,1,2}} & \vdots & {w_{N_x,N_y,2}} \\ \vdots \\ {w_{1,1,N_z}} & \vdots & {w_{N_x,N_y,N_z}}\\
  \end{bmatrix} 
\end{equation*}

The factorization with direction splitting solver is cheaper than the generation of the right-hand sides. We solve three one-dimensional problems with multiple right-hand sides. The cost of the generation of the right-hand sides is high, but it can be reduced around one order of magnitude by switching to the piece-wise constant basis functions.

We compare the standard RHS generation code

{\tt
\begin{tabbing}
xx\=xx\=xx\=xx\=xx\=xx\=xx\=xx\=xx\=xx\kill
1  \> {\bf for} nex=1,$N_x$ //loop through elements along $x$ \\
2  \> {\bf for} ney=1,$N_y$ //loop through elements along $y$ \\
3  \> {\bf for} nez=1,$N_z$ //loop through elements along $z$ \\
4  \> \> {\bf for} ibx=1,p+1 //loop through p+1 B-splines along $x$\\
5  \> \> {\bf for} iby=1,p+1 //loop through p+1 B-splines along  $y$\\
6  \> \> {\bf for} ibz=1,p+1 //loop through p+1 B-splines along  $z$\\
7  \>  \> \> i = f(nex,ibx) //global index of B-spline along $x$ \\
8  \>  \> \> j = f(ney,iby) //global index of B-spline along $y$ \\
9  \>  \> \> k = f(nez,ibz) //global index of B-spline along $z$ \\
10  \>  \> \> irow = g(nex,ibx,ney,iby,nez,ibz) // global row index \\
11  \> \> \> {\bf for} qx=1,nqx //quadrature point along $x$ \\
12  \> \> \> {\bf for} qy=1,nqy //quadrature point along $y$ \\
13  \> \> \> {\bf for} qz=1,nqz //quadrature point along $z$ \\
  \> \> \> \>   // aggregate RHS \\ 
14  \> \> \> \> L(irow)+= $weight*RHS(qx,qy,qz,B^x_{i,p}(qx)B^y_{j,p}(qy)B^z_{k,p}(qz))$ \\
\end{tabbing} }

with the one where the test functions are set to piece-wise constants

{\tt
\begin{tabbing}
xx\=xx\=xx\=xx\=xx\=xx\=xx\=xx\=xx\=xx\kill
1  \> {\bf for} nex=1,$N_x$ //loop through elements along $x$ \\
2  \> {\bf for} ney=1,$N_y$ //loop through elements along $y$ \\
3  \> {\bf for} nez=1,$N_z$ //loop through elements along $z$ \\
4  \> \> {\bf for} ibx=1,p+1 //loop through p+1 B-splines along $x$\\
5  \> \> {\bf for} iby=1,p+1 //loop through p+1 B-splines along  $y$\\
6  \> \> {\bf for} ibz=1,p+1 //loop through p+1 B-splines along  $z$\\
7  \>  \> \> irow = g(nex,ibx,ney,iby,nez,ibz) // global row index \\
8  \> \> \> {\bf for} qx=1,nqx/2 //quadrature point along $x$ \\
9  \> \> \> {\bf for} qy=1,nqy/2 //quadrature point along $y$ \\
10  \> \> \> {\bf for} qz=1,nqz/2 //quadrature point along $z$ \\
  \> \> \> \>   // aggregate RHS \\ 
11  \> \> \> \> L(irow)+= $weight*RHS(qx,qy,qz,1.0)$ \\
\end{tabbing} }

We implement the isogeometric L2 projection using quadratic $C^1$ B-splines, and piece-wise constant test functions. We report in Table 1 and Figure \ref{figure2} the cost of generation of the right-hand-sides, and the cost of factorization. 
We execute the code on a single core of a Linux cluster node with 2.4 GHz Intel Xeon CPU E5-2509. We conclude that switching to piece-wise constant test functions reduces the cost almost one order of magnitude, using the slowest traditional integration with Gaussian quadrature. The further speedup can be possibly obtained by incorporating faster quadrature \cite{Barton} and parallel solvers \cite{CAI}.

{\tiny\begin{table*}[htp]
 \resizebox{0.9\textwidth}{!}{
 \begin{subtable}[h]{0.3\textwidth}
\centering
\begin{tabular}{lll }
\hline
\multicolumn{3}{c}{ quadratic B-splines C1 } \\
$n_x=n_y=n_z$ & \#NRDOF & time[s] \\
\hline
2 &64 & 0.0013 \\
4 &216 & 0.0085 \\
8 &1,000 & 0.065 \\
16 &5,832 & 0.53 \\
32 & 39,304 & 3.96 \\
64 & 287,496& 31.24 \\
128 &2,197,000 & 250.42 \\
256 &17,173,512 &  2004.00\\
\hline
\end{tabular}
\end{subtable}
\hspace{2cm}
\begin{subtable}[h]{0.3\textwidth}
\centering
\begin{tabular}{|lll }
\hline
\multicolumn{3}{c}{ piece-wise constants } \\
$n_x=n_y=n_z$ & \#NRDOF & time[s] \\
\hline
2 &64 &  0.0005 \\
4 &216 &  0.0022 \\
8 &1,000 &  0.0143\ \\
16 &5,832 &  0.123\\
32 &39,304 &  0.64 \\
64 &287,496 &  4.90 \\
128 &2,197,000 & 39.17 \\
256 & 17,173,512 & 338.36 \\
\hline
\end{tabular}
\end{subtable}
\hspace{2cm}
 \begin{subtable}[h]{0.3\textwidth}
\centering
\begin{tabular}{|lll }
\hline
\multicolumn{3}{c}{ factorization time } \\
$n_x=n_y=n_z$ & \#NRDOF & time[s] \\
\hline
2 & 64 & 0.0005 \\
4 & 216 & 0.00009 \\
8 & 1,000 & 0.004  \\
16 & 5,832 & 0.028 \\
32 &39,304  & 0.21 \\
64 &287,496  &  1.66 \\
128 &2,197,000 & 13.32 \\
256 &17,173,512 &  106.00 \\ 
\hline
\end{tabular}
\end{subtable}}
\caption{Fortran 90 implementation of 3D isogeometric L2 projection on a cluster node. Generation time for test functions set to either quadratic B-splines with$C^1$ continuity, or piece-wise constants. Factorization time (does not depend on the generation method in case of direct solver). \#NRDOF denotes the number of degrees of freedom, $n_x,n_y,n_z$ denotes the number of elements along $x,y,z$ axes.}
\label{table1}
\end{table*}}

We also consider the improvement from the application of the piece-wise constant test functions, when we use higher-order B-splines, for quadratics, cubics, and quartics, over the larger mesh. We report the times in Table \ref{table1a}.

{\tiny\begin{table*}[htp]
 \resizebox{0.9\textwidth}{!}{
 \begin{subtable}[h]{0.3\textwidth}
\centering
\begin{tabular}{lll }
\hline
\multicolumn{3}{c}{ $C^{p-1}$ B-splines } \\
$p$ & \#NRDOF & time[s] \\
\hline
2 & 17,173,512 & 2,004 \\
3 & 17,373,979 & 10,571 \\
4 & 17,576,000 & 38,902 \\
\hline
\end{tabular}
\end{subtable}
\hspace{2cm}
\begin{subtable}[h]{0.3\textwidth}
\centering
\begin{tabular}{|lll }
\hline
\multicolumn{3}{c}{ piece-wise constants } \\
$p$ & \#NRDOF & time[s] \\
\hline
2 & 17,173,512 & 338 \\ 
3 & 17,373,979 & 667 \\
4 & 17,576,000 & 1,243 \\
\hline
\end{tabular}
\end{subtable}
\hspace{2cm}
 \begin{subtable}[h]{0.3\textwidth}
\centering
\begin{tabular}{|lll }
\hline
\multicolumn{3}{c}{ factorization time } \\
$p$ & \#NRDOF & time[s] \\
\hline
2 & 17,173,512 & 106 \\ 
3 & 17,373,979 & 234\\
4 & 17,576,000 & 420 \\
\hline
\end{tabular}
\end{subtable}}
\caption{Fortran 90 implementation of 3D isogeometric L2 projection on a cluster node. Generation time for test functions set to either $C^{p-1}$ B-splines, or piece-wise constants. Factorization time (does not depend on the generation method in case of direct solver). Mesh size is fixed for $256\times 256 \times 256$ elements, and $p$ denotes the B-splines order.}
\label{table1a}
\end{table*}}

\begin{figure}
\includegraphics[scale=0.5]{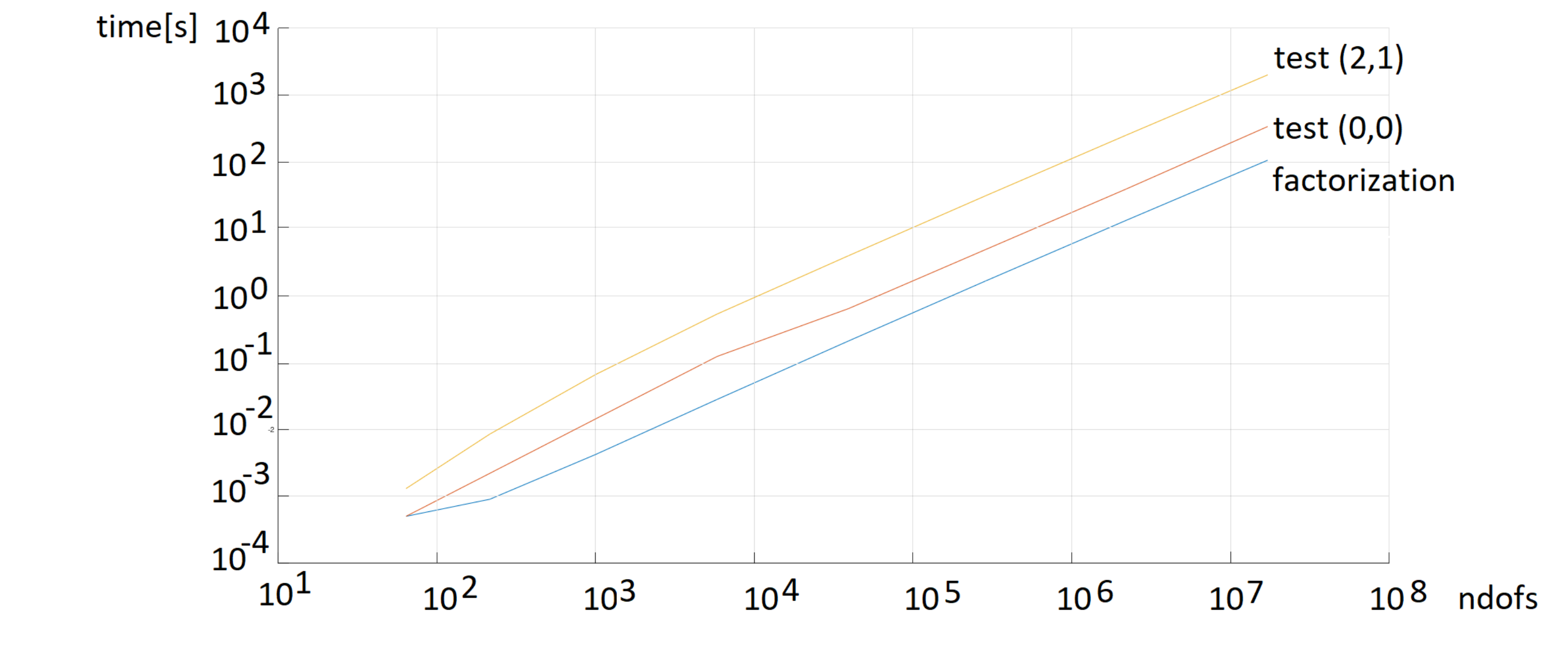}
\caption{Fortran 90 implementation of 3D isogeometric L2 projection on a cluster node, with quadratic B-splines (denoted by (2,1)) and piece-wise constant B-splines (denoted by (0,0)). Factorization by alternating directions solver.}
\label{figure2}
\end{figure}

\subsection{Explicit dynamics}

We focus now on the time-dependent problems solved with an explicit method.
The governing equation in the strong form is given by
\begin{eqnarray}\frac{\partial u}{\partial t} - {\cal F} u = {\cal RHS} \end{eqnarray}

The strong form is transformed into a weak one by taking the $L^2$ scalar product with test functions $v \in H^1\left(\Omega\right)$,
and the Euler integration scheme is utilized with respect to time
\begin{eqnarray}\left(v,u_{t+1}\right)_{L^2} = 
\left(v, u_t+ Dt {\cal F} u_t + Dt{\cal RHS}\right)_{L^2} \end{eqnarray}

The system has an identical structure as the one considered in the projection problem, and the "elimination" of test functions can be applied here as well, speeding up the integration at every time step.

{\footnotesize \begin{equation*}
	\begin{bmatrix}
    \int {\cal I}_{1}{\cal I}_{1}{\cal I}_{1}{\left(B^x_{1,p}B^y_{1,p}B^z_{1,p}\right)} & \cdots &  \int {\cal I}_{1}{\cal I}_{1}{\cal I}_{1}{\left(B^x_{N_x,p}B^y_{N_y,p}B^z_{N_z,p}\right)} \\
    \int {\cal I}_{2}{\cal I}_{1}{\cal I}_{1}{\left(B^x_{1,p}B^y_{1,p}B^z_{1,p}\right)}& \cdots &  \int {\cal I}_{2}{\cal I}_{1}{\cal I}_{1}{\left(B^x_{N_x,p}B^y_{N_y,p}B^z_{N_z,p}\right)} \\   
    \vdots & \vdots & \vdots \\
        \int {\cal I}_{N_x}{\cal I}_{N_y}{\cal I}_{N_z}{\left(B^x_{1,p}B^y_{1,p}B^z_{1,p}\right)} & \cdots &  \int {\cal I}_{N_x}{\cal I}_{N_y}{\cal I}_{N_z}{\left(B^x_{N_x,p}B^y_{N_y,p}B^z_{N_z,p}\right)} \\
  \end{bmatrix}
  \begin{bmatrix}
  {u_{1,1,1}} \\ {u_{2,1,1}} \\ \vdots \\ {u_{N_x,N_y,N_z}}\\
  \end{bmatrix} 
\end{equation*}
\begin{equation*} 
=    \begin{bmatrix}
\int {\cal I}_{1}{\cal I}_{1}{\cal I}_{1} \left( u_t+ Dt {\cal F} u_t + {\cal RHS}(x,y,z)  \right)dxdydz \\
\int {\cal I}_{2}{\cal I}_{1}{\cal I}_{1} \left( u_t+ Dt {\cal F} u_t + {\cal RHS}(x,y,z)  \right)dxdydz \\
\vdots \\
\int {\cal I}_{N_x}{\cal I}_{N_y}{\cal I}_{1} \left(  u_t+ Dt {\cal F} u_t + {\cal RHS}(x,y,z)\right)   dxdydz \\
  \end{bmatrix}
\end{equation*}}

We can employ the alternating directions solver in every time step.
We factorize the L2 projection matrix once, using three one-dimensional systems with multiple right-hand sides, and then we perform forward and backward substitutions for each new right-hand side. Each time step of the explicit dynamics simulation generates the right-hand-side, like in the isogeometric L2 projection problem. Thus, to get the cost of the explicit dynamics simulation, we multiply the times from Table 1 by the number of time steps. The further speedup can be obtained by using parallel explicit dynamics solvers \cite{CAI} and a fast integration scheme \cite{Barton}, reducing the number of quadrature points for the trial functions.

\subsection{Laplace problem with mixed boundary conditions} 

We consider a Laplace problem with Dirichlet and Neumann boundary conditions,  
\begin{equation}
-\Delta u = {\cal RHS},
\end{equation}
where $\Omega=\left(0,1\right)^2$, with boundary conditions
\begin{eqnarray}
u=0 \textrm{ on } \Gamma_D \\
 \frac{\partial u}{\partial n} = {\cal G} \textrm{ on } \Gamma_N
\end{eqnarray}

The weak variational formulation is obtained by taking the
$L^2$-scalar product with functions $v\in
H^1_{\Gamma_D}\left(\Omega\right)=\{v\in
H^1\left(\Omega\right):v|_{\Gamma_D}=0\}$, integrating by parts, and including the Neumann boundary conditions:
\begin{eqnarray}
\label{eq:weak1a}
\textrm{Find } u\in V = H^1_{\Gamma_D}\left(\Omega\right) \textrm{ such that} \\
b\left(v,u\right)=l\left(v\right), \forall v\in V,
\end{eqnarray}
where
\begin{eqnarray}
\label{eq:weak1B}
b\left(v,u\right)=\int_{\Omega} \nabla v \cdot \nabla u d{\bf x},
\mbox{ and} \\
l\left(v\right)=\int_\Omega v {\cal RHS} d{\bf x} +\int_{\Gamma_N} v {\cal G} dS
\label{eq:weak4b}
\end{eqnarray}

It is possible to show that the Galerkin problem is well-possed.

Now, we discretize with $C^1$ B-splines, so our solution lives in a space that is a sub-set of $H^2$, so we can integrate back by parts on a discrete level

\begin{eqnarray}
-\int_{\Omega} v_h \Delta u_h d{\bf x} +{ \int_{\Gamma_N} v_h \frac{\partial u_h}{\partial n} dS} = \int_\Omega v_h {\cal RHS} d{\bf x} + {\int_{\Gamma_N} v_h {\cal G} dS} 
\label{eq:weak4c}
\end{eqnarray}

The system in a discrete form reads

{\footnotesize \begin{equation*}
-	\begin{bmatrix}
    \int{\left(B^x_{1,p}B^y_{1,p}\right)}{\Delta\left(B^x_{1,p}B^y_{1,p}\right)} & \cdots &  \int{\left(B^x_{1,p}B^y_{1,p}\right)}{\Delta\left(B^x_{N_x,p}B^y_{N_y,p}\right)} \\
    \int{\left(B^x_{2,p}B^y_{1,p}\right)}{\Delta\left(B^x_{1,p}B^y_{1,p}\right)} & \cdots &  \int{\left(B^x_{2,p}B^y_{1,p}\right)}{\Delta\left(B^x_{N_x,p}B^y_{N_y,p}\right)} \\   
    \vdots & \vdots & \vdots \\
        \int{\left(B^x_{N_x,p}B^y_{N_y,p}\right)}{\Delta\left(B^x_{1,p}B^y_{1,p}\right)} & \cdots &  \int{\left(B^x_{N_x,p}B^y_{N_y,p}\right)}{\Delta\left(B^x_{N_x,p}B^y_{N_y,p}\right)} \\
  \end{bmatrix}
+ 
\end{equation*}
\begin{equation*}
	\begin{bmatrix}
    \int_{\Gamma_N}{\left(B^x_{1,p}B^y_{1,p}\right)}{\frac{\partial \left(B^x_{1,p}B^y_{1,p}\right)}{\partial n}} & \cdots &  \int{\left(B^x_{1,p}B^y_{1,p}\right)}{\frac{\partial \left(B^x_{N_x,p}B^y_{N_y,p}\right)}{\partial n}} \\
    \int_{\Gamma_N}{\left(B^x_{2,p}B^y_{1,p}\right)}{\frac{\partial \left(B^x_{1,p}B^y_{1,p}\right)}{\partial n}} & \cdots &  \int_{\Gamma_N}{\left(B^x_{2,p}B^y_{1,p}\right)}{{\frac{\partial\left(B^x_{N_x,p}B^y_{N_y,p}\right)}{\partial n}}} \\   
    \vdots & \vdots & \vdots \\
        \int_{\Gamma_N}{\left(B^x_{N_x,p}B^y_{N_y,p}\right)}{\frac{\partial\left(B^x_{1,p}B^y_{1,p}\right)}{\partial n}} &  \cdots &  \int_{\Gamma_N}{\left(B^x_{N_x,p}B^y_{N_y,p}\right)}{\frac{\partial\left(B^x_{N_x,p}B^y_{N_y,p}\right)}{\partial n}} \\
  \end{bmatrix}
  \begin{bmatrix}
  {u_{1,1}} \\ {u_{2,1}} \\ \vdots \\ {u_{N_x,N_y}}\\
  \end{bmatrix} 
\end{equation*}
\begin{equation*} 
=    \begin{bmatrix}
\int {\cal RHS}(x,y) {B^x_{1,p}(x)B^y_{1,p}(y)}  \\
\int {\cal RHS}(x,y) {B^x_{2,p}(x)B^y_{1,p}(y)}  \\
\vdots \\
\int {\cal RHS}(x,y) {B^x_{N_x,p}(x)B^y_{N_y,p}(y)} \\
  \end{bmatrix}
+ \begin{bmatrix}
\int_{\Gamma_N} {B^x_{1,p}(x)B^y_{1,p}(y)} {\cal G}(x,y) \\
\int_{\Gamma_N} {B^x_{2,p}(x)B^y_{1,p}(y)} {\cal G}(x,y) \\
\vdots \\
\int_{\Gamma_N} {B^x_{N_x,p}(x)B^y_{N_y,p}(y)}{\cal G}(x,y) \\
  \end{bmatrix}
\end{equation*}}

Now, we move to the piece-wise constant test functions

\begin{equation*}
-	\begin{bmatrix}
  \int {\cal I}_{1}{\cal I}_{2}{\Delta\left(B^x_{1,p}B^y_{1,p}\right)}  & \cdots &  \int {\cal I}_{1}{\cal I}_{2}{\Delta\left(B^x_{N_x,p}B^y_{N_y,p}\right)} \\  
    \vdots & \vdots &  \vdots \\
  \int {\cal I}_{i}{\cal I}_{j}{\Delta\left(B^x_{1,p}B^y_{1,p}\right)}  & \cdots &  \int {\cal I}_{i}{\cal I}_{j}{\Delta\left(B^x_{N_x,p}B^y_{N_y,p}\right)} \\  
    \vdots & \vdots &  \vdots \\
        \int {\cal I}_{N_x}{\cal I}_{N_y}\Delta\left(B^x_{1,p}B^y_{1,p}\right) & \cdots &  \int {\cal I}_{N_x}{\cal I}_{N_y}{\Delta\left(B^x_{N_x,p}B^y_{N_y,p}\right)} \\
  \end{bmatrix}
\end{equation*}
\begin{equation*}
	+\begin{bmatrix}
    \int {\cal I}_{1}{\cal I}_{1}\cap\Gamma_N{\frac{\partial \left(B^x_{1,p}B^y_{1,p}\right)}{\partial n}} & \cdots &  \int {\cal I}_{1}{\cal I}_{1}\cap\Gamma_N{\frac{\partial \left(B^x_{N_x,p}B^y_{N_y,p}\right)}{\partial n}} \\
    \vdots & \vdots &  \vdots \\
    \int {\cal I}_{i}{\cal I}_{j}\cap\Gamma_N{\frac{\partial \left(B^x_{1,p}B^y_{1,p}\right)}{\partial n}} & \cdots & \int {\cal I}_{i}{\cal I}_{j}\cap\Gamma_N{{\frac{\partial\left(B^x_{N_x,p}B^y_{N_y,p}\right)}{\partial n}}} \\   
    \vdots & \vdots & \vdots \\
        \int {\cal I}_{N_x}{\cal I}_{N_y}\cap\Gamma_N{\frac{\partial\left(B^x_{1,p}B^y_{1,p}\right)}{\partial n}} &  \cdots &  \int {\cal I}_{N_x}{\cal I}_{N_y}\cap\Gamma_N{\frac{\partial\left(B^x_{N_x,p}B^y_{N_y,p}\right)}{\partial n}} \\
  \end{bmatrix}
  \begin{bmatrix}
  {u_{1,1}} \\ \vdots \\ {u_{i,j}} \\ \vdots \\ {u_{N_x,N_y}}\\
  \end{bmatrix} 
\end{equation*}
\begin{equation*} 
=    \begin{bmatrix}
\int {\cal I}_{1}{\cal I}_{2}{\cal RHS}(x,y) \\
\vdots \\
\int {\cal I}_{i}{\cal I}_{j}{\cal RHS}(x,y) \\
\vdots \\
\int {\cal I}_{N_x}{\cal I}_{N_y}{\cal RHS}(x,y)  \\
  \end{bmatrix}
+ \begin{bmatrix}
\int {\cal I}_{1}{\cal I}_{1}\cap\Gamma_N {\cal G} \\
\vdots \\
\int {\cal I}_{i}{\cal I}_{j}\cap \Gamma_N  {\cal G} \\
\vdots \\
\int {\cal I}_{N_x}{\cal I}_{N_y}\Gamma_N {\cal G} \\
  \end{bmatrix}
\end{equation*}
The zero Dirichlet boundary conditions can be enforced by setting corresponding rows to 0, diagonals to 1, and right-hand-sides to 0.

We compare the standard aggregation code

{\tt
\begin{tabbing}
xx\=xx\=xx\=xx\=xx\=xx\=xx\=xx\=xx\=xx\kill
1  \> {\bf for} nex=1,$N_x$ //loop through elements along $x$ \\
2  \> {\bf for} ney=1,$N_y$ //loop through elements along $y$ \\
3  \> \> {\bf for} ibx1=1,p+1 //loop through p+1 B-splines along $x$\\
4  \> \> {\bf for} iby1=1,p+1 //loop through p+1 B-splines along  $y$\\
5  \>  \> \> i = f(nex,ibx1) //global index of B-spline along $x$ \\
6  \>  \> \> j = f(ney,iby1) //global index of B-spline along $y$ \\
7  \>  \> \> irow = g(nex,ibx1,ney,iby1) // global row index \\
8  \> \> \> {\bf for} qx=1,nqx //quadrature point along $x$ \\
9  \> \> \> {\bf for} qy=1,nqy //quadrature point along $y$ \\
  \> \> \> \>   // aggregate RHS \\ 
10  \> \> \> \> L(irow)+= $weight*RHS(qx,qy,B^x_{i,p}(qx)B^y_{j,p}(qy))$ \\
11  \> \> \> \> {\bf for} ibx2=1,p+1 //loop through p+1 element B-splines along $x$\\
12  \> \> \> \> {\bf for} iby2=1,p+1 //loop through p+1 element B-splines along  $y$\\
13  \>  \> \> \> \> k = f(nex,ibx2) //global index of B-spline along $x$ \\
14  \>  \> \> \> \> l = f(ney,iby2) //global index of B-spline along  $y$ \\
15 \>  \> \> \> \> icol = g(nex,ibx2,ney,iby2) // global column index \\
16 \> \> \> \> \>  {\bf for} rx=1,nqx //quadrature point along $x$ \\
17 \> \> \> \> \> {\bf for} ry=1,nqy //quadrature point along $y$ \\
 \> \> \> \> \>  \> // aggregate LHS \\
18  \> \> \> \> \> \> M(irow,icol)+= weight*$A(B^x_{i,p}(qx)B^y_{j,p}(qy),B^x_{k,p}(rx)B^y_{l,p}(ry))$  \\
\end{tabbing} }

with the one where the test functions are set to piece-wise constants

{\tt
\begin{tabbing}
xx\=xx\=xx\=xx\=xx\=xx\=xx\=xx\=xx\=xx\kill
1  \> {\bf for} nex=1,$N_x$ //loop through elements along $x$ \\
2  \> {\bf for} ney=1,$N_y$ //loop through elements along $y$ \\
3  \> \> {\bf for} ibx1=1,p+1 //loop through p+1 B-splines along $x$\\
4  \> \> {\bf for} iby1=1,p+1 //loop through p+1 B-splines along  $y$\\
5  \>  \> \> irow = g(nex,ibx1,ney,iby1) // global row index \\
6  \> \> \> {\bf for} qx=1,nqx/2 //quadrature point along $x$ \\
7  \> \> \> {\bf for} qy=1,nqy/2 //quadrature point along $y$ \\
  \> \> \> \>   // aggregate RHS \\ 
8  \> \> \> \> l(irow)+= $weight*RHS(qx,qy,1.0)$ \\
9  \> \> \> \> {\bf for} ibx2=1,p+1 //loop through p+1 piece-wise constant along $x$\\
10  \> \> \> \> {\bf for} iby2=1,p+1 //loop through p+1 piece-wise constant along  $y$\\
11 \>  \> \> \> \> icol = g(nex,ibx2,ney,iby2) //global column index \\
 \> \> \> \> \>  // aggregate LHS \\
12  \> \> \> \> \> M(irow,icol)+= $weight*A(B^x_{i,p}(qx)B^y_{j,p}(qy),1.0)$  \\
\end{tabbing} }

Namely, we verify the execution times using the MATLAB implementation executed on a laptop. 
The comparison is presented in Table \ref{table2} and Figure \ref{figure2}. Further reduction of the execution time can be obtained by using fast quadrature \cite{Barton} or parallel integration \cite{CPC}.

{\tiny\begin{table*}[htp]
 \resizebox{0.9\textwidth}{!}{
 \begin{subtable}[h]{0.3\textwidth}
\centering
\begin{tabular}{lll }
\hline
\multicolumn{3}{c}{ quadratic B-splines C1 } \\
$n_x=n_y$ & \#NRDOF & time[s] \\
\hline
4 & 36 & 2 \\
8 & 100 &  9 \\
16 & 324 &  35 \\
32 & 1,156 &  130 \\
64 & 4,356 &  521 \\
128 & 16,900 &  2100\\
256 & 66,564 & 8204\\
\hline
\end{tabular}
\end{subtable}
\hspace{2cm}
\begin{subtable}[h]{0.3\textwidth}
\centering
\begin{tabular}{|lll}
\hline
\multicolumn{3}{c}{ piece-wise constants } \\
$n_x=n_y$ & \#NRDOF & time[s] \\
\hline
4  & 36 & 0.1 \\
8 & 100 & 0.5\\
16 & 324 & 2 \\
32 & 1,156 & 9 \\
64 & 4,356 & 34 \\
128 & 16,900 &  131 \\
256 & 66,564 &  523 \\
\hline
\end{tabular}
\end{subtable}
\hspace{2cm}able
 \begin{subtable}[h]{0.3\textwidth}
\centering
\begin{tabular}{|lll }
\hline
\multicolumn{3}{c}{ factorization time } \\
$n_x=n_y$ & \#NRDOF & time[s] \\
\hline
4  & 36 & 0.0009 \\
8 & 100 &   0.004 \\
16 & 324 &  0.02 \\
32 & 1,156 &  0.15 \\
64 & 4,356 & 1.19 \\
128 & 16,900 &  10.04 \\
256 & 66,564 & 70.53 \\
\hline
\end{tabular}
\end{subtable}}
\caption{MATLAB implementation of the 2D Laplace problem on a laptop. Generation time for test functions set to either quadratic B-splines with $C^1$ continuity, or piece-wise constants. Factorization time (does not depend on the generation method in case of direct solver). \#NRDOF denotes the number of degrees of freedom, $n_x,n_y$ denotes the number of elements along $x,y$ axes.}
\label{table2}
\end{table*}}

\begin{figure}
\includegraphics[scale=0.5]{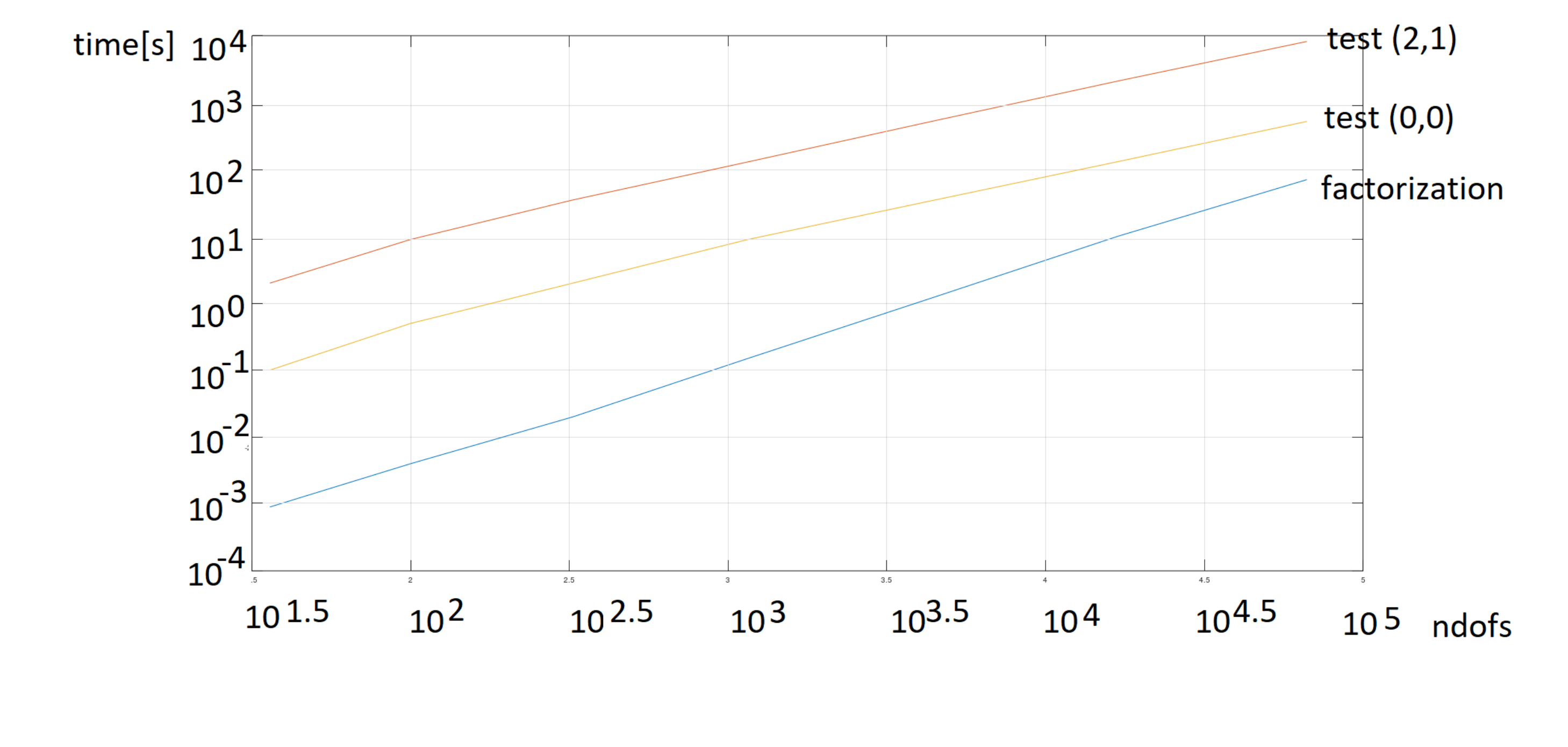}
\caption{MATLAB implementation of 2D Laplace problem on a laptop, with quadratic B-splines (denoted by (2,1)) and piece-wise constant B-splines (denoted by (0,0)). Factorization by MATLAB "backslash" solver.}
\label{figure2}
\end{figure}

\subsection{Isogeometric L2 projection of a bitmap} 

Finally we consider the isogeometric L2 projection of a bitmap. We decompose the bitmap into three RGB tables with [0,255] values denoting the contributions from the red, green and blue colors. We solve the three projection problems, and we combine the results to get the colors.

\begin{equation}
u = {\cal BITMAP},
\end{equation}

We present the resulting bitmaps, obtained by executing our method with piece-wise constant test functions and quadratic $C^1$ trial B-splines. We also present in Figure \ref{figure:curves} the comparison of our method with the isogeometric L2 projection with quadratic $C^1$ B-splines for trial and test.

\begin{figure}
\centering
\includegraphics[scale=0.3]{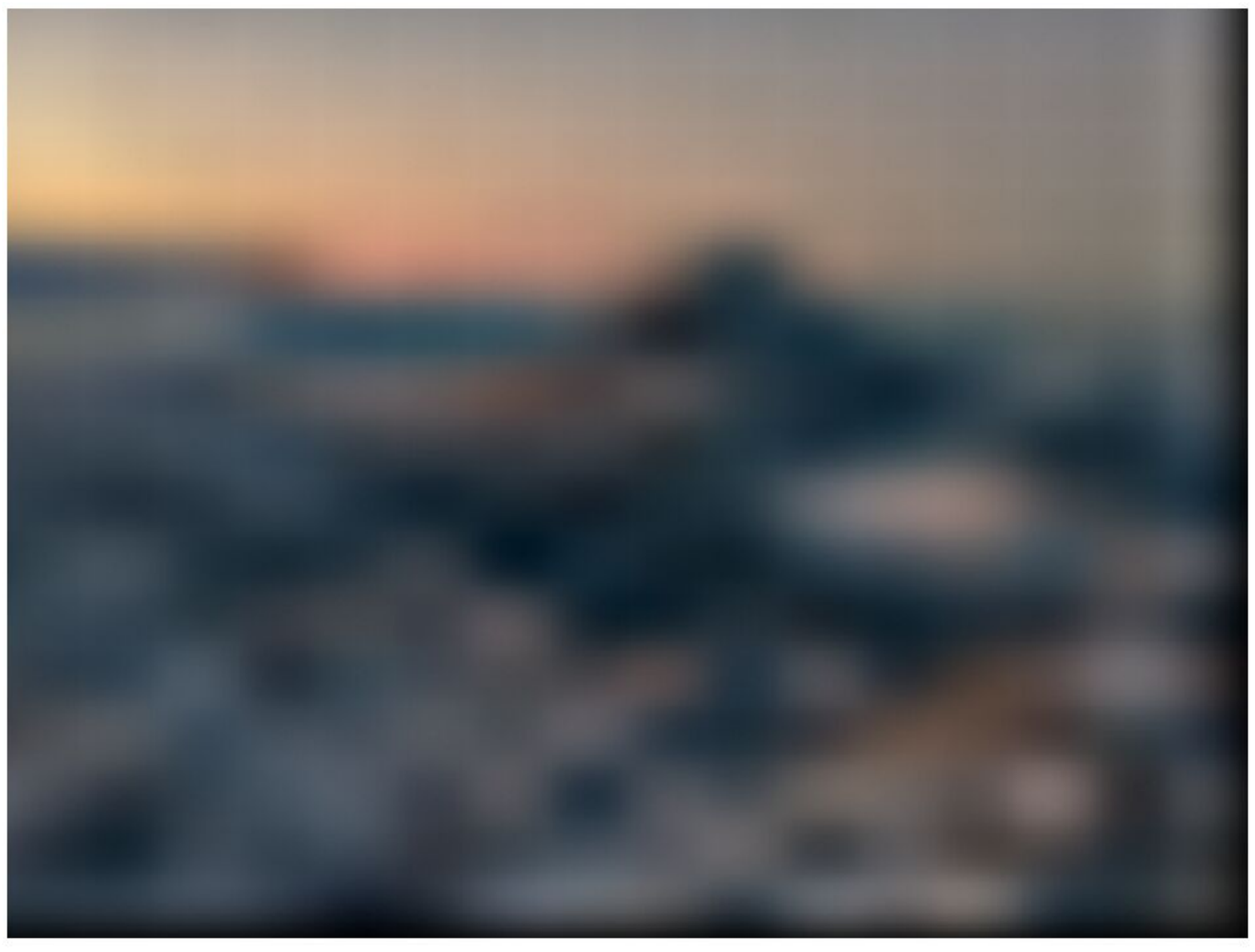}
\includegraphics[scale=0.3]{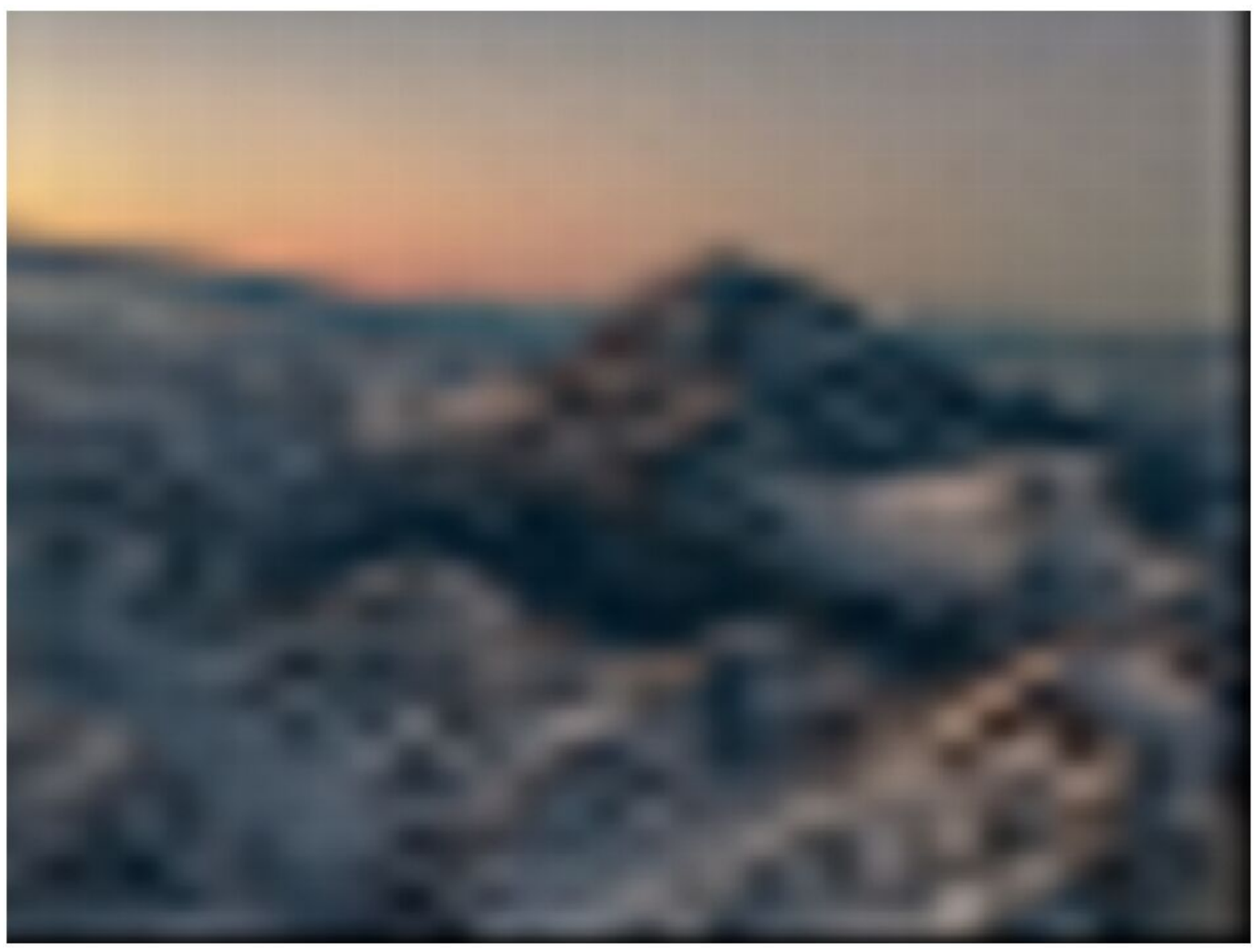}
\includegraphics[scale=0.3]{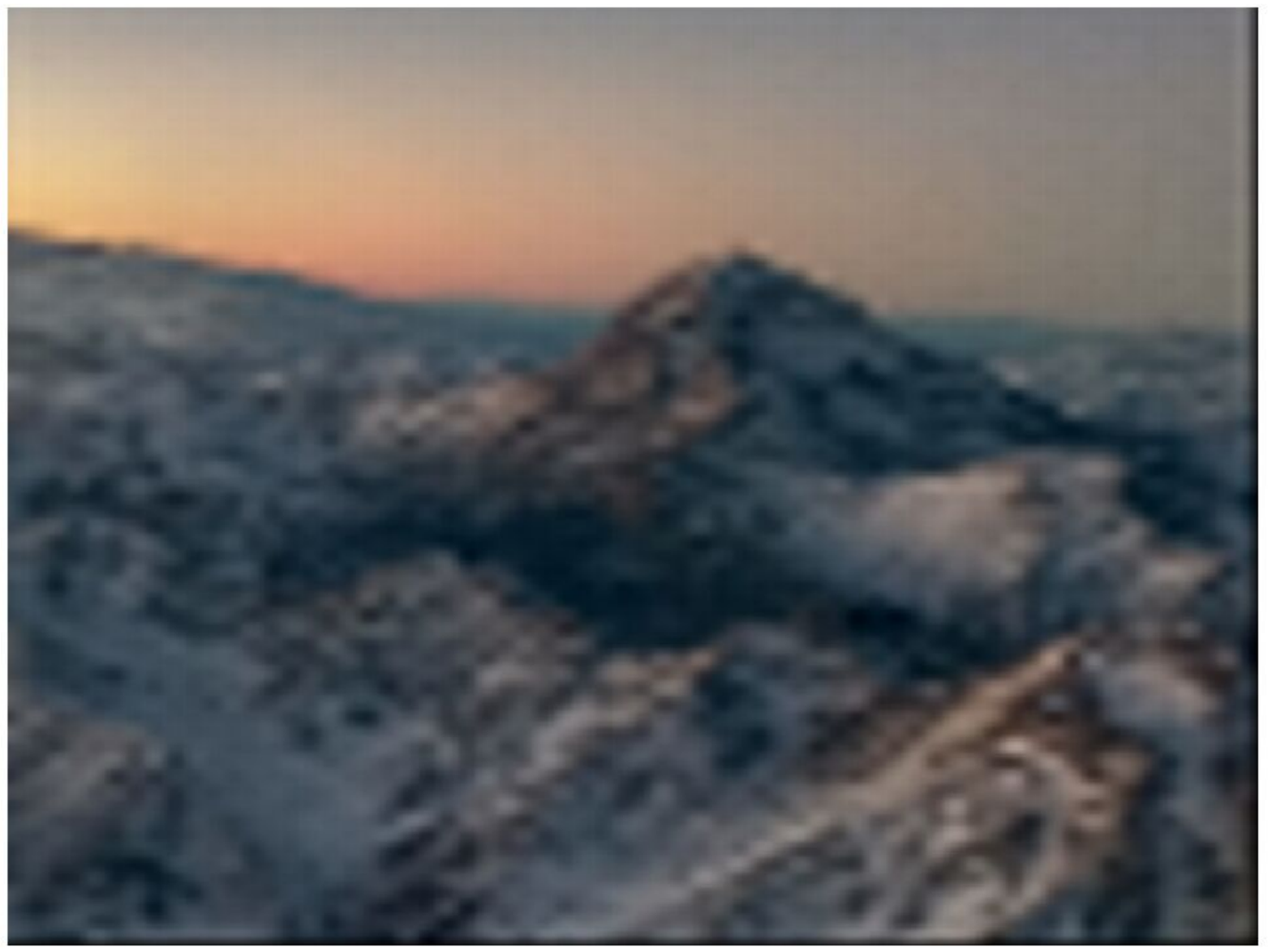}
\includegraphics[scale=0.3]{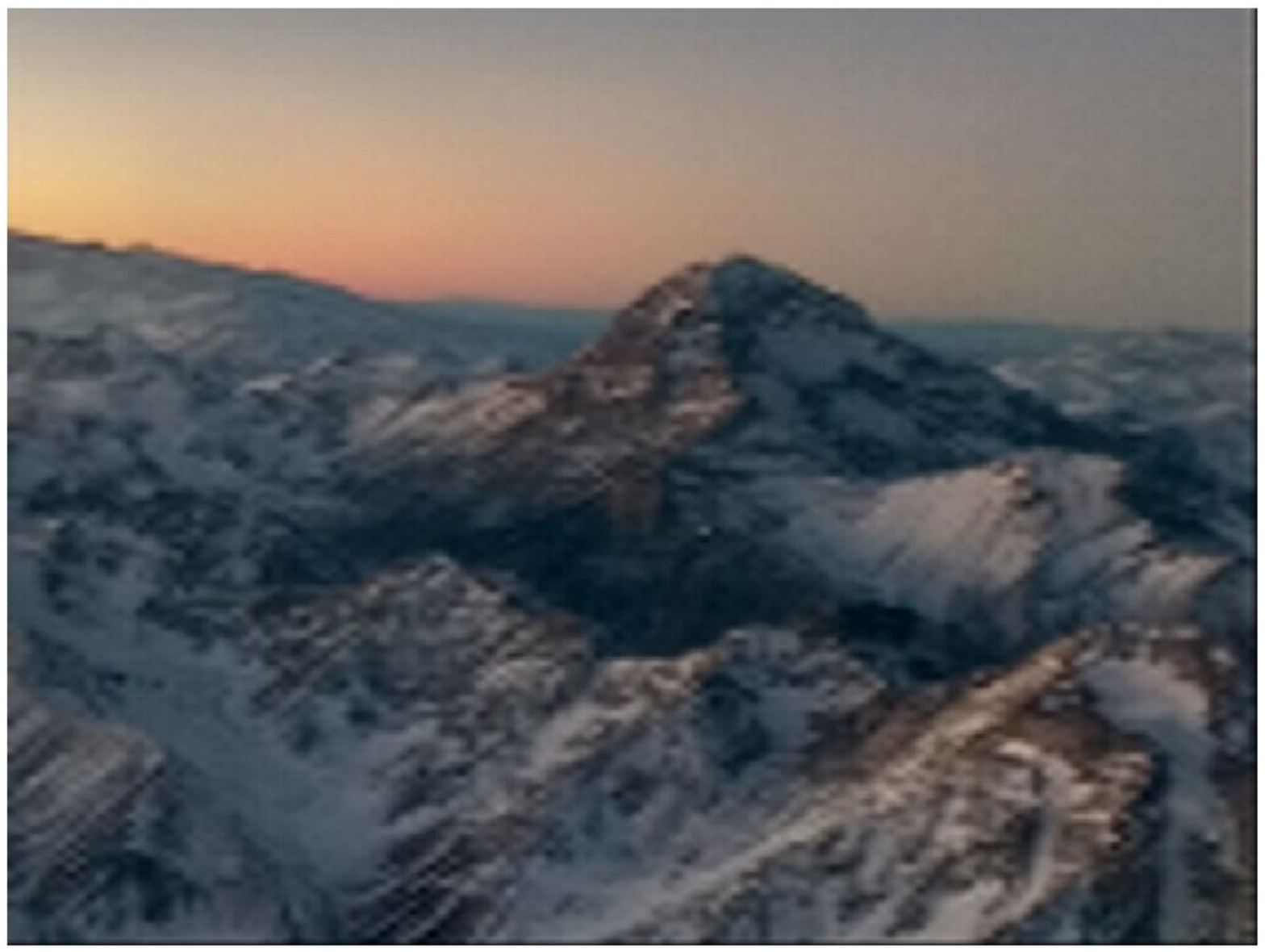}
\caption{MATLAB implementation of the isogeometric L2 projection problem for a bitmap, using $16\times16, 32\times 32, 64\times 64$ and $128\times 128$ meshes with quadratic B-splines and piece-wise constant polynomials, span over particular elements of trial space.}
\label{figureAnia}
\end{figure}

\begin{figure}
\centering
\includegraphics[scale=0.6]{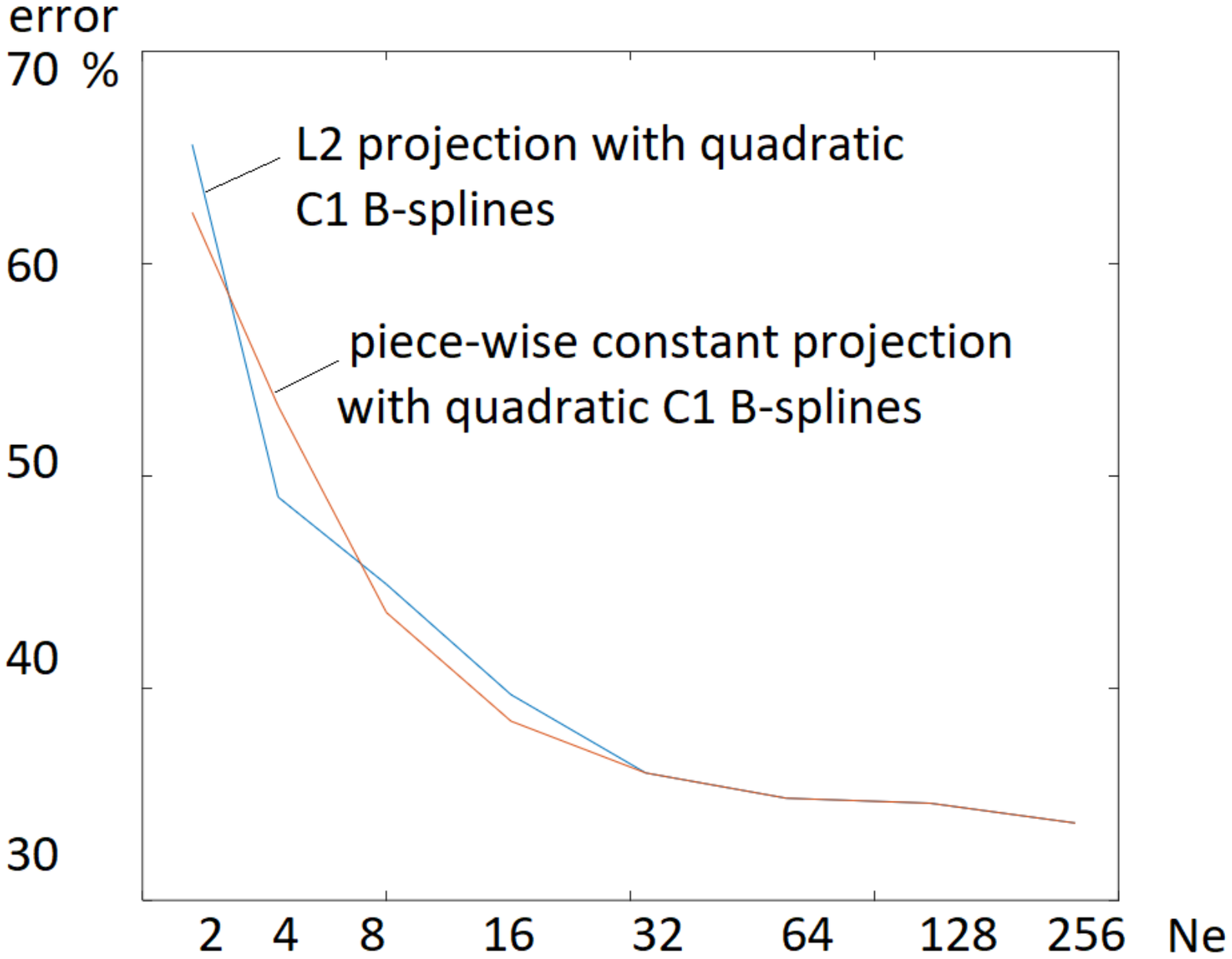}
\caption{Comparison of convergence of the isogeometric L2 projection of the bitmaps with quadratic $C^1$ B-splines for trial and test, and with our method with piece-wise constant test functions and quadratic $C^1$ B-splines for approximation. The exact error is measured in $L^2$ norm.}
\label{figure:curves}
\end{figure}

\section{Conclusions}

We have shown in this paper, that solving a PDE with Galerkin method with $H^2$ approximation of $C^1$ basis functions, can be transformed into testing the PDE with piece-wise constant test functions. The resulting problem is of the Petrov-Galerkin kind, with different trial and test spaces.
This has the following consequences. First, we can eliminate the test functions from the linear systems of equations, by making them piece-wise constants. 
Second, the numerical integration cost will be reduced since we do not need to integrate the test functions right-away. Third, this method is PDE independent, but we cannot integrate by parts since the derivatives of B-splines do not fulfill the partition of unity property at a given quadrature point. 
However,  when we use higher continuity, e.g., $C^1$ discretizations, the system of equations integrated by parts is equivalent to the system not integrated by parts (the entries in the matrices are indeed equal).
Our method is of Petrov-Galerkin kind, where we discretize with higher continuity basis functions preserving the partition of unity property, and test with piece-wise constant functions.
Fourth, the method does not depend on the selected quadrature. Fifth, the method does not depend on the shape of the domain. Sixth, the method is dimension independent, and it can be used in space-time formulations as well. Finally, the method can be used to speed up IGA time-dependent simulations.

Summing up, the test functions in IGA can be set to a piece-wise constant. The test functions define the span of the collocation points. The points are combined with the weights as prescribed by the quadrature for the integration. The collocations are computed at the quadrature points, and they are combined over the span of the test functions. The same logic applies to any basis functions that are globally $C^1$, and they preserve the partition of unity property. 

Future work may include the mathematical analysis of this new projection method. We will also check how "removing" of test functions from IGA discretizations influences the convergence of iterative solvers \cite{IT}. We will also check how it does affect space-time formulations \cite{ST}. This method can also be combined with some fast integration techniques (in the sense that we only reduce the order of the integrated function, so our method does not exclude further speedup by using faster quadrature). We also plan to investigate how this method can be incorporated with some stabilization methods \cite{S1,S2,S3,S4}.

\section*{Acknowledgments}
This work is supported by National Science Centre, Poland grant no. 2017/26/M/ ST1/ 00281. 
I would like to thank prof. David Pardo and dr Marcin \L{}o\'{s} for discussion on the limitations of the method.

\end{document}